\documentclass{amsart}
\usepackage{amsmath}
\usepackage{amsthm}
\usepackage{amsfonts}
\usepackage{amssymb}
\usepackage[T1]{fontenc}
\usepackage{hyperref}
\usepackage[capitalise]{cleveref}
\usepackage{esint}
\usepackage{mathrsfs, mathtools, bm}
\usepackage{xcolor}
\usepackage{comment}
\usepackage{tikz}
\theoremstyle{plain}
\usetikzlibrary{calc,shapes.geometric,positioning,arrows.meta}
\usepackage{fullpage}
\usepackage{lmodern}
\usepackage{microtype}
\newtheorem{theorem}{Theorem}[section]
\newtheorem{defn}[theorem]{Definition}
\newtheorem{proposition}[theorem]{Proposition}
\newtheorem{corollary}[theorem]{Corollary}
\newtheorem{lemma}[theorem]{Lemma}

\newtheorem{remark}[theorem]{Remark}
\newcommand{\R}{\mathbb{R}}

\newcommand{\Sp}{\mathbb{S}}

\title{Affine boosted Gauss curvature flow}
\author{Ben Andrews}
\address{Mathematical Sciences Institute, Australian National University, Canberra ACT 2601, Australia}
\email{ben.andrews@anu.edu.au}
\author{Devesh Rajpal}
\address{Mathematical Sciences Institute, Australian National University, Canberra ACT 2601, Australia}
\email{devesh.rajpal@anu.edu.au}
\date{}
\subjclass[2020]{35K55, 53A05}
\keywords{Gauss curvature flow, Isotropic measure, Nonlinear parabolic flow} 
\begin{document}

\begin{abstract}
    We prove convexity in the affine orbit of entropies associated to Gauss curvature flow for sub-affine powers. This allows the construction of an affine boosted Gauss curvature flow, which keeps the evolving hypersurface at the affine minimizer of the Gaussian entropy. The modified flow has bounded isoperimetric ratio, and after rescaling it converges along a sequence of times to self-similar solutions of the Gauss curvature flow. The construction also proves existence of new self-similar solutions to Gauss curvature flow in the sub-affine regime.
\end{abstract}
\maketitle
\section{Introduction}

Firey considered an isotropic model of strictly convex tumbling stones on a beach where the speed of a point on the surface is equal to the Gauss curvature in the inward normal direction \cite{firey1974shapes}. Assuming regularity, uniqueness and antipodal symmetry he proved the flow converges to a round sphere. Since Firey, a lot of work has been done on Gauss curvature flow including flow by powers of Gauss curvature. A one-parameter family of hypersurfaces \( X:M^{n} \times [0,T) \to \R^{n+1} \) is said to evolve by \( \alpha \)-Gauss curvature flow (\( \alpha \)-GCF) if it satisfies 
\[ \frac{\partial}{\partial t}X(p,t) = - K^{\alpha}(p,t)\nu(p,t). \]

Tso proved regularity, short-time existence and convergence to a point of \( \alpha \)-GCF \cite{tso1985deforming}. Firey's conjecture was proved by the first author for surfaces \cite{andrews1999gauss}. 
The affine exponent \( \alpha = 1 / (n+2) \) is the affine normal flow for which the flow is invariant under special linear group action. The first author established convergence of the rescaled flow to ellipsoids \cite{andrews1996contraction} in this case. For \( \alpha = 1/n \), Chow \cite{chow1985deforming} proved convergence to sphere. Other authors have established convergence to a self-similar solution for the range \( \alpha > 1/(n+2) \) \cite{andrews2000motion,guan2017entropy,andrews2016flow,chow1985deforming}. It was proved in \cite{brendle2017asymptotic} that spheres are the only self-similar solutions in the super-affine range, which settled a conjecture of Firey in all dimensions and exponent range.

In the sub-affine regime \( \alpha < 1/(n+2) \), for generic initial data the reaction term dominates the diffusion term for the rescaled flow leading to unbounded isoperimetric ratio along the flow. This was established for origin symmetric curves in \cite{andrews2002non} and extended to non-symmetric curves for a dense open set of initial data in \cite{andrews2003classification}. Nonetheless, it is an interesting question to ask what are the non-spherical solutions of the fully non-linear elliptic PDE 
\[ \left< X,\nu \right> = K^{\alpha} \]
which characterize self-similar solutions of the flow contracting to the origin. The solution of the above equation can also be viewed as solution of the \( L_{p} \)-Minkowski problem for \( p = 1 - 1 / \alpha \).
The one dimensional case has been fully classified by the first author \cite{andrews2003classification}. The solutions are \( k \)-symmetric curves for \( \alpha < 1/(k^{2}-1)\). In higher dimensions, the first author established the existence of self-similar solutions with irreducible symmetry group \( G \subset SO(n+1) \) \cite{andrews2000motion}. The goal of this paper is to establish existence of new self-similar solutions in the sub-affine regime with weaker symmetry assumptions. In particular, we prove the following theorem. 
\begin{theorem}\label{maintheorem}
    Let \( G \subset O(n+1) \) be a closed subgroup that acts non-transitively on \( \Sp^{n} \). Then there is an integer \( d \ge 3 \) such that for every \( \alpha \in \big( 0, \frac{1}{(d-1)(n+d)} \big) \) there exists a non-spherical, smooth, compact, uniformly convex self-similar solution to the \( \alpha \)-GCF whose symmetry group contains \( G \).
\end{theorem}
In particular, this gives non-spherical self-similar solutions with axial symmetry for every \( \alpha \in (0, \frac{1}{2(n+3)}) \) (\cref{axial}), and solutions whose symmetry group contains a given finite group, such as a dihedral group. This resolves a conjecture from \cite{andrews2000motion} in the affirmative.

We introduce an affine corrected version of sub-affine GCF to prove this. The existence of monotone entropy for a large class of geometric parabolic flows was established in \cite{andrews1997monotone}. The entropy concerned in this paper is the translation invariant Gaussian entropy given by 
\[ Z_{1}(\mathcal{M}) = \left( \frac{V(\mathcal{M})}{V(B^{n+1})} \right)^{n(\alpha -1)} \left( \fint_{\Sp^{n}}K^{\alpha - 1}d \mu \right)^{n+1} \]
which is non-decreasing under \( \alpha \)-GCF.

Special affine positions have been explored in convex geometry under various setting. Petty proved the existence of surface area minimizing position \cite{petty1961surface}. Giannopoulos-Milman generalized it to quermassintegrals \cite{giannopoulos1999isotropic,giannopoulos2000extremal}. Lutwak-Yang-Zhang further explored the \( L_{p} \)-surface area minimizing positions \cite{Lutwak_Yang_Zhang_2005}. We observe that for GCF, the entropies exhibit strict convexity in the pure affine orbit in the sub-affine regime. This implies existence of a minimizing affine position. Specifically, for any compact locally uniformly convex hypersurface, there exists a \( T_{0} \in \operatorname{SL}(n+1) \) such that \[ Z_{1}(T_{0}\mathcal{M}) \le Z_{1}(T \mathcal{M}) \] for any \( T \in \operatorname{SL}(n+1) \) with equality iff \( T \) is related to \( T_{0} \) by an orthogonal transformation. The position \( T_{0}\mathcal{M} \), for \( T_{0} \) symmetric,  will be called \( Z_{1} \)-isotropic position and enjoys special geometric properties; in particular, its entropy is bounded by a dimensional constant. The convexity of entropy in affine orbit may be of independent interest. This is covered in \cref{convexityaffine,geometryatZ1}.

Utilizing this special affine position, we modify the \( \alpha \)-GCF by a global affine term,
\[ \frac{\partial X}{\partial t} = -K^{\alpha}\nu + \Lambda(t)X, \qquad \Lambda(t) \in \mathfrak{sl}(n+1)/\mathfrak{so}(n+1), \]
where \( \Lambda(t) \) is chosen so that the hypersurface stays at the \( Z_{1} \)-isotropic position (see \cref{sec:ABGCF} for the precise definition). We call this flow the Affine Boosted Gauss Curvature Flow (ABGCF). Our second main result is the following.
\begin{theorem}\label{introconvergence}
    Let \( 0 < \alpha < \frac{1}{n+2} \) and let \( \mathcal{M}_{0} = \partial\Omega_{0} \) be a smooth, closed, locally uniformly convex hypersurface in \( \R^{n+1} \) in \( Z_{1} \)-isotropic position. Then the ABGCF starting at \( \mathcal{M}_{0} \) has a unique smooth solution on a maximal time interval \( [0,T) \) with \( T < \infty \), and
    \begin{enumerate}
        \item the isoperimetric ratio \( \rho_{+}(\mathcal{M}_{t})/\rho_{-}(\mathcal{M}_{t}) \) is bounded by a constant depending only on \( n \), \( \alpha \) and \( Z_{1}(\mathcal{M}_{0}) \);
        \item \( \operatorname{diam}(\Omega_{t}) \to 0 \) as \( t \to T \);
        \item there exist a curve \( p(t) \in \overline{\Omega}_{t} \) and times \( t_{i} \to T \) such that the hypersurfaces \[\big( V(B^{n+1})/V(\Omega_{t_{i}}) \big)^{\frac{1}{n+1}}\big( \mathcal{M}_{t_{i}} - p(t_{i}) \big) \] converge in \( C^{\infty} \) to a smooth, uniformly convex hypersurface \( \mathcal{M}_{\infty} \) which is a self-similar solution of the \( \alpha \)-GCF, is in \( Z_{1} \)-isotropic position, and satisfies \( Z_{1}(\mathcal{M}_{\infty}) \ge Z_{1}(\mathcal{M}_{0}) \).
    \end{enumerate}
\end{theorem}
The proof of a theorem like \cref{introconvergence} in the literature is assembled by proving curvature estimates. The new ingredient is the global term \( \Lambda(t) \) which needs to be controlled via \textit{a priori} estimates. We prove an estimate of \( |\Lambda| \) against the time derivative of entropy and other geometric parameters in \cref{xilemma}, which is sufficient for curvature estimates of ABGCF in \cref{sec:lowerspeed,sec:upperspeed}.
\Cref{maintheorem} follows by applying \cref{introconvergence} to small \( G \)-invariant perturbations of the sphere, chosen so that the initial entropy exceeds that of the sphere; this is done in \cref{construction}. 

\subsection*{AI Disclosure.} Most of the work in this paper was done by early 2026 without AI. Later, some parts of the proofs in Section 8 and 9 were polished using AI.

\section{Notation and preliminary results}

In this section we set up some notation related to convex bodies in Euclidean space and discuss established results which will be needed later.

Let \( \Omega \subset \R^{n+1} \) be a convex body, that is a convex set with non-empty interior. Associated to each convex body \( \Omega \) is a support function \( \sigma_{\Omega} : \Sp^{n} \to \R \) defined by 
\[ \sigma_{\Omega}(z) = \sup_{p \in \Omega} \left< p, z \right>. \]
For \( \Omega \) locally uniformly convex and with smooth boundary, the Gauss map \( z \to \nu(X(z)) \) has an inverse \( X: \Sp^{n} \to \R^{n+1} \) in terms of support function as the following expression 
\[ X(z) = \sigma(z)z + \nabla\sigma(z). \]
We can also extract the second fundamental form from the support function. Define an operator \( A[f]_{ij} = \nabla_{i}\nabla_{j}f + g_{ij}f \), where \( \nabla \) is the standard covariant derivative and \( g \) is the standard metric on \( \Sp^{n} \). It is known that the eigenvalues of \( A[\sigma] \) with respect to \( g \) are the principal radii \( r_{1} = \kappa_{1}^{-1}, \dots, r_{n} = \kappa_{n}^{-1} \) of the hypersurface at point \( X(z) \). Thus, the Gauss curvature \( K = \kappa_{1} \cdots \kappa_{n} \) can be written as \( K = \det g / \det A[\sigma] \). Often we will fix a local orthonormal frame on \( \Sp^{n} \) and use the expression \( K(X(z)) = \det A[\sigma]^{-1}(z) \). By divergence formula, the volume is given by 
\[ V(\Omega) = \frac{1}{n+1}\int_{\partial \Omega} \langle X, \nu \rangle d \mathcal{H}^{n} = \frac{1}{n+1} \int_{\Sp^{n}} \frac{\sigma(z)}{K(z)}d \mu = \frac{1}{n+1} \int_{\Sp^{n}}\sigma \det A[\sigma] d \mu.\]

The Minkowski addition of convex body \( \Omega_{1} \) and \( \Omega_{2} \) is given by \( \Omega_{1}+\Omega_{2} = \{p_{1} + p_{2} : p_{1} \in \Omega_{1}, p_{2} \in \Omega_{2}\} \), which is another convex body. It can be verified that \( \sigma_{\Omega_{1}+ \Omega_{2}} = \sigma_{\Omega_{1}} + \sigma_{\Omega_{2}} \). This implies that volume of convex bodies is a homogenous polynomial of degree \( n+1 \) under Minkowski addition. Let \( \Omega_{0}, \dots, \Omega_{n} \) be \( (n+1) \) bounded convex bodies, define mixed volume \( V(\Omega_{0}, \ldots, \Omega_{n}) \) by the formula 
\[ V(t_{0}\Omega_{0} + \cdots + t_{n}\Omega_{n}) = \sum_{0 \le i_{0},\dots, i_{n} \le n}t_{i_{0}} \ldots t_{i_{n}}V(\Omega_{i_{0}}, \dots, \Omega_{i_{n}}).  \]
Using the volume formula, we can also write the mixed volume as 
\[ V(\Omega_{0},\dots, \Omega_{n}) = \frac{1}{n+1}\int_{\Sp^{n}}\sigma_{0}Q[\sigma_{1},\dots, \sigma_{n}]d \mu, \]
where \( Q[\sigma_{1},\dots, \sigma_{n}] \) is the mixed determinant operator given by 
\[ Q[\sigma_{1}, \dots, \sigma_{n}] = \frac{1}{n!} \sum_{\tau, \omega \in S_{n}} (-1)^{|\tau|}(-1)^{|\omega|}A[\sigma_{1}]_{\tau(1)}^{\omega(1)} \cdots A[\sigma_{n}]_{\tau(n)}^{\omega(n)},  \]
where \( S_{n} \) denote the permutation group with \( n \) elements and the operator \( A \) has raised indices using the metric \( g \) on \( \Sp^{n} \). We can recover the determinant by setting \( \sigma_{1} = \cdots = \sigma_{n} = \sigma \), yielding \( Q[\sigma,\dots,\sigma] = \det A[\sigma]. \)
We will often write \( Q_{n} = Q[\sigma,\dots, \sigma] \) for brevity when it is clear from the context. For fixed bodies \( \Omega_{2},\dots, \Omega_{n} \) the mixed determinant's derivative can also be viewed as an elliptic operator, defined by the formula 
\[ Q[f, \sigma_{2}, \dots, \sigma_{n}] = \dot{Q}^{ij}A[f]_{ij} = \dot{Q}^{ij}(\nabla_{i}\nabla_{j}f + g_{ij}f) \]
for any smooth \( f \). The operator \( \dot{Q}^{ij} \) is divergence-free. When \( \sigma_{2} = \cdots = \sigma_{n} = \sigma \), we have
\begin{equation}
    \dot{Q}^{ij} = \frac{1}{n}Q_{n}[\sigma](A[\sigma]^{-1})^{ij}, \label{Qdot}
\end{equation}
since \( \frac{d}{ds}\det A[\sigma + sf]\big|_{s=0} = nQ[f,\sigma,\dots,\sigma] = Q_{n}(A^{-1})^{ij}A[f]_{ij} \). The reader can refer to \cite[Section 18.7]{andrews2022extrinsic} for more details. A very important theorem on volume of convex bodies is the Brunn-Minkowski theorem. 
\begin{theorem}[Brunn-Minkowski Theorem]\label{BM}
    Let \( \Omega_{0}, \Omega_{1} \) be convex bodies in \( \R^{n+1} \), then the function \( V^{\frac{1}{n+1}}  \) is concave with respect to Minkowski addition, that is for any \( t \in [0,1] \), 
    \[ V((1-t)\Omega_{0} + t\Omega_{1})^{\frac{1}{n+1}} \ge (1-t)V(\Omega_{0})^{\frac{1}{n+1}} + tV(\Omega_{1})^{\frac{1}{n+1}},  \]
    with equality if and only if \( \Omega_{0},\Omega_{1} \) are related by homothety and translation.
    The concavity is equivalent to the following inequality for any convex bodies \( \Omega_{0},\Omega_{1} \), 
    \[ V(\Omega_{0},\Omega_{1},\dots, \Omega_{1})^{2} \ge V(\Omega_{0},\Omega_{0}, \Omega_{1}, \dots,\Omega_{1})V(\Omega_{1},\dots,\Omega_{1}).  \]
    If \( \Omega_{1} \) has smooth boundary and is uniformly convex, equality holds if and only if \( \Omega_{0} \) and \( \Omega_{1} \) are related by homothety and translation. Assuming the convex bodies have smooth boundary, expanding into support functions we get the following Poincar\'e type inequality on sphere, for any smooth function \( f \) on \( \Sp^{n} \), 
    \[ \int_{\Sp^{n}} f^{2} Q[1,\sigma_{1},\dots, \sigma_{1}]d \mu \le \int_{\Sp^{n}}\nabla_{i}f \nabla_{j}f \dot{Q}^{ij}d \mu + \frac{(\int_{\Sp^{n}} f Q[\sigma_{1},\dots, \sigma_{1}]d \mu)^{2}}{\int_{\Sp^{n}}\sigma_{1}Q[\sigma_{1},\dots, \sigma_{1}]d \mu}, \]
    with equality if and only if \( f(z) = \left< p, z \right> + c \sigma_{1}(z) \) for some \( p \in \R^{n+1} \) and \( c \in \R \).
\end{theorem} 

Affine optimization problems on convex bodies have been studied before in \cite{petty1961surface,giannopoulos1999isotropic,giannopoulos2000extremal}. The critical points of this action lead to the concept of isotropic measures, which in our setting can be defined as follows.
\begin{defn}
    A non-negative finite Borel measure \( \mu \) on \( \Sp^{n} \) is isotropic if 
    \[  \int_{\Sp^{n}}|x \cdot z|^{2}d \mu(z) = \frac{|\mu|}{n+1} \]
    for all \( x \in \Sp^{n} \) where \( |\mu| = \int_{\Sp^{n}}d \mu \). Equivalently, if for any linear map \( T :\R^{n+1} \to \R^{n+1}\), 
    \[ \int_{\Sp^{n}}\left< Tz,z \right>d \mu(z) =  \frac{\operatorname{tr}(T)|\mu|}{n+1}. \]
\end{defn}
Examples of isotropic measures include the spherical Lebesgue measure and cross measure given by \( \mu = \sum_{i=1}^{n+1}\delta_{e_{i}} + \sum_{i=1}^{n+1}\delta_{-e_{i}}   \) where \( \{e_{i}\} \) is an orthonormal frame of \( \R^{n+1} \). 

\section{Entropy under affine orbit} \label{convexityaffine}

A classic result of Petty states that a convex body has minimal surface area in its affine orbit if and only if, its area measure is isotropic \cite{petty1961surface}. The result has been generalized to $L_p$ surface area measures as well by Lutwak, Yang and Zhang in their celebrated result on existence of $L_p$ John ellipsoid \cite{Lutwak_Yang_Zhang_2005}. The isotropy of the self-similar solutions of \( \alpha \)-GCF was also established by Hu-Ivaki \cite{hu2025stability} and used to prove uniqueness in the range \( \alpha \in (1/2, 1) \). We will show that there is a range of integrals associated to a smooth convex body which are similarly minimized under affine orbit if and only if, the associated measure is isotropic.  Examples include Firey and Gaussian entropy in the sub-affine range which are monotone under \( \alpha \)-GCF. Following \cite{andrews1997monotone} we define a range of integrals associated to the \( \alpha \)-GCF.

Let \( \beta >0 \), and let \( \Omega \) contain the origin in its interior if \( \beta \neq 1 \). Define the scale invariant quantity \( Z_{\beta} \) as
\begin{align*}
    Z_{\beta}(\sigma) & = \left( \frac{V(\Omega)}{V(B^{n+1})}\right)^{\beta -1 + n(\alpha \beta -1)} \left(  \fint_{\Sp^{n}} \left( \frac{1}{\sigma Q_{n}[\sigma]^{\alpha}}\right)^{\beta} \sigma Q_{n}[\sigma] d \mu\right)^{n+1} \\
    & = \left( \fint_{\Sp^{n}}\sigma Q_{n}[\sigma]d \mu\right)^{\beta -1 + n(\alpha \beta -1)}\left( \fint_{\Sp^{n}}\sigma^{1-\beta}Q_{n}[\sigma]^{1-\alpha \beta}d \mu \right)^{n+1}.
\end{align*}
\begin{defn}
    We say that \( \Omega \) is in \( Z_{\beta} \)-isotropic position if the measure \( \sigma^{1-\beta}Q_{n}[\sigma]^{1-\alpha\beta}\,d\mu \) is isotropic.
\end{defn}
In order to understand the behavior of \( Z_{\beta} \) under affine orbit, we derive the formulas for change in support function \( \sigma \) and curvature function \( Q_{n}[\sigma] \) under affine transformation. For general convex bodies, the reader can look at \cite[Lemma 12]{schutt2004surface} or \cite[Lemma 7.12]{lutwak1990centroid}. For \( C^{2} \) convex bodies, we adapt a proof from \cite[Section 16.1]{andrews2022extrinsic}.

\begin{lemma}
    Let \( \Omega \subset \R^{n+1} \) be a convex body with smooth boundary \( \mathcal{M} = \partial \Omega \) and \( L \in \operatorname{GL}^{+}(n+1) \). The support function and curvature function of \( L\Omega \) are given by 
    \[ \sigma_{L\Omega}(z) = \sigma_{\Omega}(\left< L^{T}z \right>)|L^{T}z|, \qquad Q_{n}[\sigma_{L\Omega}](z) = (\det L)^{2}Q_{n}[\sigma_{\Omega}](\left< L^{T}z \right>)|L^{T}z|^{-n-2},  \]
    where \( L^{T} \) denotes the transpose and \( \left< x \right> = x/|x| \).
    \label{affineQ}
\end{lemma}
\begin{proof}
    Let \( \{\partial_{i}\} \subset T_{X(z)}\mathcal{M} \) be a normal coordinate system of \( \mathcal{M} \) at \( X(z) \) with normal \( z \) such that \( \partial_{1} \wedge \cdots \wedge \partial_{n} \wedge z = \omega \) is the volume form of \( \R^{n+1} \). The vectors \( \{L\partial_{i}\} \) form a basis of \( T_{LX(z)}L\mathcal{M} \). Let \( z^{L} \) denote the normal of \( L\mathcal{M} \) at \( LX(z) \), then \( z^{L} \) satisfies the equality \( 0 = \left< z^{L}, L \partial_{i} \right> = \left< L^{T}z^{L}, \partial_{i} \right> \) for all \( i \). Thus, \( z^{L} = L^{-T}z / |L^{-T}z| = \left< L^{-T}z \right>\) where \( L^{-T} \) denote the inverse transpose. We also have the identity \( |L^{-T}z| = |L^{T}z^{L}|^{-1} \).

    The support function of \( L\Omega \) at \( z^{L} \) is given by 
    \[ \sigma_{L\Omega}(z^{L}) = \left< LX(z), z^{L} \right> = \left< X(z), L^{T}z^{L} \right> = \sigma_{\Omega}(z)|L^{-T}z|^{-1} = \sigma_{\Omega}(z)|L^{T}z^{L}|. \] 
    The second fundamental form of \( L\mathcal{M} \) at \( LX(z) \) is given by 
    \[ h_{ij}^{L}(z^{L}) = -\left< \frac{\partial}{\partial x^{i}}(L\partial_{j}), z^{L} \right> = -\left< L\frac{\partial}{\partial x^{i}}\partial_{j}, z^{L} \right> = -\left< \frac{\partial}{\partial x^{i}}\partial_{j}, L^{T}z^{L} \right> = h_{ij}(z)|L^{-T}z|^{-1}. \]
    The metric of \( L\mathcal{M} \) at \( LX(z) \) is given by 
    \[ g_{ij}^{L}(z^{L}) = \left< L\partial_{i}, L\partial_{j} \right>. \] 
    Note that \( L \partial_{1} \wedge \cdots \wedge L \partial_{n} \wedge z^{L} = \sqrt{\det g_{ij}^{L}} \, \omega \) and \( L \partial_{1} \wedge \cdots \wedge L \partial_{n} \wedge Lz = \det L \, \omega \). The component of \( Lz \) in the direction of \( z^{L} \) is given by \( \left< Lz, z^{L} \right> = |L^{-T}z|^{-1} \). It follows that 
    \( \sqrt{\det g_{ij}^{L}}|L^{-T}z|^{-1} = \det L \), or equivalently \( \det g_{ij}^{L} = (\det L)^{2}|L^{-T}z|^{2} \). Thus, the Gauss curvature of \( L\mathcal{M} \) at \( LX(z) \) is given by 
    \[ K^{L}(z^{L}) = \frac{\det h_{ij}^{L}}{\det g_{ij}^{L}} = \frac{\det h_{ij}|L^{-T}z|^{-n}}{\det g_{ij}(\det L)^{2}|L^{-T}z|^{2}} = K(z)(\det L)^{-2}|L^{-T}z|^{-n-2}. \]
    This also implies \( Q_{n}[\sigma_{L \Omega}](z^{L}) = (\det L)^{2}Q_{n}[\sigma_{\Omega}](z) \)\( |L^{-T}z|^{n+2} = (\det L)^{2}Q_{n}[\sigma_{\Omega}](z)|L^{T}z^{L}|^{-n-2} \).
\end{proof}
We are now ready to prove the main proposition of this section.

\begin{proposition}
    Let \( \alpha <\frac{1}{n+2} \) and \( \beta > 0\). Let \( \Omega \subset \R^{n+1} \) be a fixed locally uniformly convex body with the origin in its interior if \( \beta \neq 1 \). Then there exists a unique positive-definite symmetric matrix \( T \in \operatorname{SL}(n+1) \) such that \( T \Omega \) is in isotropic \( Z_{\beta} \) position. Moreover,
    \[ Z_{\beta}(T \Omega) \le Z_{\beta}(S \Omega) \qquad \forall S \in \operatorname{SL}(n+1),\]
    with equality if and only if \( S = OT \) for some \( O \in \operatorname{SO}(n+1) \). \label{isotropicZ}
\end{proposition}
\begin{proof}
    Since volume isn't changed by \( \operatorname{SL}(n+1) \) transformations, we can ignore the volume factor of \( Z_{\beta} \) for now. Let \( \tilde{Z} = \fint_{\Sp^{n}}\sigma^{1-\beta}Q_{n}^{1- \alpha \beta} \). The affine group \( \operatorname{SL}(n+1) \) acts on the space of convex bodies and the Lie algebra \( \mathfrak{sl}(n+1) \) generates one-parameter group of convex bodies. Let \( A \in \mathfrak{sl}(n+1) / \mathfrak{so}(n+1) \) be a symmetric traceless matrix. Let \( \Omega_{s} = e^{sA}\Omega \). Denote the support function of \( \Omega_{s} \) by \( \sigma_{s} \). Using \( \sigma_{s}(z) = \sigma_{0}(e^{sA}z / |e^{sA}z|)|e^{sA}z|  \) we can derive the differential change as
    \begin{equation*}
        D_{A}\sigma(z) = \frac{d}{ds}\bigg|_{s=0}\sigma_{s}(z) = \sigma \left< Az,z \right> + \left< \nabla \sigma, (Az)^{T}\right>.
    \end{equation*}
    Let \( \xi(z) = \left< Az,z \right> \), then \( \nabla \xi = 2(Az)^{T} \). It follows that 
    \begin{equation}
        D_{A}\sigma(z) = \sigma \xi + \frac{1}{2}\left< \nabla \sigma, \nabla \xi \right>. \label{DAsigma}
    \end{equation}
    Similarly, from \( Q_{n}[\sigma_{s}](z) = Q_{n}[\sigma_{0}](e^{sA}z/|e^{sA}z|) |e^{sA}z|^{-n-2}\) we get 
    \begin{equation}
      D_{A}Q_{n}[\sigma](z) = \frac{d}{ds}\bigg|_{s=0}Q_{n}[\sigma_{s}] = -(n+2)Q_{n}\xi + \frac{1}{2}\left< \nabla Q_{n}, \nabla \xi \right>. \label{DAQ}
    \end{equation}
    
    Using the affine transformation formulas it is easy to see that 
    \[ \tilde{Z}(\sigma_{s}) = \fint_{\Sp^{n}}\sigma_{0}(\left< e^{sA}z \right>)^{1- \beta}Q_{n}[\sigma_{0}](\left< e^{sA}z \right>)^{1-\alpha \beta}|e^{sA}z|^{1-\beta - (n+2)(1-\alpha \beta)}d \mu(z). \]
    Let \( y = \left< e^{sA}z\right> = \frac{e^{sA}z}{|e^{sA}z|} \), then the Jacobian of the transformation is \( d \mu(y) = |e^{sA}z|^{-n-1}d \mu(z) = |e^{-sA}y|^{n+1}d \mu(z) \). It follows that 
    \begin{equation}
        \tilde{Z}(\sigma_{s}) = \fint_{\Sp^{n}}\sigma_{0}^{1-\beta}Q_{n}[\sigma_{0}]^{1-\alpha \beta}|e^{-sA}y|^{\beta(1- \alpha (n+2))}d \mu(y).\label{Zs}
    \end{equation}
    Let \( \gamma = \beta(1- \alpha(n+2)) \) and define \( f(s) = |e^{-sA}y|^{\gamma} \) for a fixed \( y \in \Sp^{n} \). The first and second derivative of \( f \) are 
    \[ f'(s) = \gamma|e^{-sA}y|^{\gamma-2}\left< e^{-sA}y,e^{-sA}(-A)y \right>, \]
    and 
    \begin{align*}
        f''(s) & = \gamma |e^{-sA}y|^{\gamma-4}\left((\gamma-2)\left< e^{-sA}y, e^{-sA}(-A)y\right>^{2} + 2 |e^{-sA}y|^{2}|e^{-sA}(-A)y|^{2} \right) \\
        & = \gamma^{2}|e^{-sA}y|^{\gamma-4} \left< e^{-sA}y,e^{-sA}Ay \right>^{2} + 2\gamma|e^{-sA}y|^{\gamma-4} \left(|e^{-sA}y|^{2}|e^{-sA}Ay|^{2} - \left< e^{-sA}y,e^{-sA}Ay \right>^{2} \right)
    \end{align*}
    Notice that the second term in the last-line is non-negative by Cauchy-Schwarz and rearranging gives \( f(s)f''(s) - f'(s)^{2} \ge 0 \), with equality only if \( y \) is an eigenvector of \( A \), which is a set of measure zero. Thus, \( \tilde{Z}(\sigma_{s}) \) is a log-convex function of \( s \). Moreover, from \cref{Zs} one has \( \tilde{Z}(\sigma_{s}) \to \infty \) as \( s \to \pm\infty \) so \( Z \) has a unique minimum in direction \( A \) over \( \operatorname{SL}(n+1) \). In order to prove properness over \( \mathfrak{sl}(n+1) / \mathfrak{so}(n+1) \), for any \( c>0 \), consider the set \( P = \{A \in \mathfrak{sl}(n+1) / \mathfrak{so}(n+1) : \tilde{Z}(e^{A}\Omega) \le c\} \). If there exists a sequence \( A_{k} \in P \) such that \( |A_{k}| \to \infty \), then we can find a subsequence \( A_{k_{j}} \) such that \( A_{k_{j}}/|A_{k_{j}}| \to A_{\infty} \) for some \( A_{\infty} \in \mathfrak{sl}(n+1) / \mathfrak{so}(n+1) \). Then, by the previous argument, \( \tilde{Z}(e^{sA_{\infty}}\Omega) \to \infty \) as \( s \to \pm\infty \). This contradicts the assumption that \( |A_{k}| \to \infty \) and \( A_{k} \in P\). Thus, the set \( P \) is bounded and closed in \( \mathfrak{sl}(n+1) / \mathfrak{so}(n+1) \).
    
    It follows that \( \tilde{Z} \) is a proper convex function in \( \operatorname{SL}(n+1) / \operatorname{SO}(n+1) \) and has a unique minimum in \( \operatorname{SL}(n+1) \) up to orthogonal transformations. If \( \Omega \) is at a minimizing position, then for any \( A \in \mathfrak{sl}(n+1) / \mathfrak{so}(n+1) \),
    \[0 = D_{A}\tilde{Z} = \frac{d}{ds}\bigg|_{s=0}\tilde{Z}(\sigma_{s}) = (- \gamma) \fint_{\Sp^{n}}\sigma^{1-\beta}Q_{n}[\sigma]^{1-\alpha \beta}\left< y,Ay \right>d \mu(y). \]
    The isotropic condition follows from setting \( A = E_{ij} + E_{ji} \) and \( A = E_{ii} - E_{jj} \) for \( i \neq j \) where \( E_{ij} \) is the elementary matrix with entry \( 1 \) at \( (i,j) \) position and zero at other positions. Conversely, a critical point of the strictly convex function \( \log\tilde{Z} \) on \( \operatorname{SL}(n+1)/\operatorname{SO}(n+1) \) is its minimum, so \( Z_{\beta} \)-isotropic positions are exactly the minimizers.
\end{proof}

\begin{remark}
    The proof is also valid for \( \alpha < 0 \); however, we will not consider expanding flows in this paper.
\end{remark}
\begin{remark}
    Substituting \( \alpha = 0 \) and \( \beta = p \) recovers the result of Lutwak-Yang-Zhang \cite{Lutwak_Yang_Zhang_2005} on minimizing \( L_{p} \) curvature over affine orbit. Also, for \( 0 < \alpha < \frac{1}{n+2} \) and \( \beta = \frac{1}{\alpha} \), the integral corresponds to the generalized dual quermassintegral of the polar body \( \Omega^{*} \) whose associated isotropic positions were considered by Bastero-Romance \cite{bastero2004positions}.
\end{remark}

\section{Affine Boosted Gauss Curvature Flow}\label{sec:ABGCF}
We will be working primarily with Gaussian entropy 
\begin{equation}
Z_{1} = \left( \frac{V(\Omega)}{V(B^{n+1})} \right)^{n(\alpha-1)}\left(\fint_{\Sp^{n}}Q_{n}[\sigma]^{1-\alpha}\right)^{n+1} \label{Z1def}
\end{equation} which is translation invariant and strictly increasing under \( \alpha \)-GCF unless \( \sigma \) is a self-similar solution \cite{andrews2022extrinsic}.
\begin{defn}
    A smooth hypersurface \( \mathcal{M} \) is said to be at \( Z_{1} \)-isotropic  position if \( dZ_{1} = Q_{n}^{1-\alpha}d \mu \) is isotropic.
\end{defn}
By \cref{isotropicZ} there always exists a \( Z_{1} \)-isotropic position in the affine orbit up to orthogonal transformation. We will use this to create a Gauss curvature flow with global terms such that the hypersurface stays at the \( Z_{1} \)-isotropic position.

\begin{defn} Let \( 0 <\alpha < 1 / (n+2) \) and \( \mathcal{M}_{0} \) a locally uniformly convex, compact hypersurface in the \( Z_{1} \)-isotropic position. The \( \alpha \)-Affine Boosted Gauss Curvature Flow (ABGCF) is a one-parameter family of hypersurfaces \( X : M^{n} \times [0,T) \to \R^{n+1} \) with \( X_{0}(M^{n}) = \mathcal{M}_{0} \) satisfying 
\begin{equation}
    \frac{\partial X}{\partial t} = - K^{\alpha}\nu + \Lambda(t)X  \label{flow1}
\end{equation}
where \( \Lambda(t) \in \mathfrak{sl}(n+1)/ \mathfrak{so}(n+1) \) is determined by the condition that \( \mathcal{M}_{t} \) is at \( Z_{1} \)-isotropic position for \( t \in [0,T) \) or equivalently \( D_{L}\tilde{Z}_{1}(\mathcal{M}_{t}) = 0 \) for all \( t \) and for all \( L \in \mathfrak{sl}(n+1) \). Identifying \( \Lambda = [\Lambda_{ab}] \) with degree two spherical harmonic \( \xi(z) = \left< \Lambda z,z  \right> = \Lambda_{ab}z^{a}z^{b} \) and using \cref{DAsigma} we can write the flow in terms of support function as 
\begin{equation}
    \frac{\partial \sigma}{\partial t} = -Q_{n}[\sigma]^{-\alpha} + D_{\Lambda}\sigma = - Q_{n}[\sigma]^{-\alpha} + \sigma \xi + \frac{1}{2}\left< \nabla \sigma, \nabla \xi \right> \label{flow2}
\end{equation}
using a time-dependent reparametrization of \( M^{n} \) which fixes the normal in time.
\end{defn}

We can find \( \xi(t) = \left< \Lambda(t)z,z \right> \) in terms of a basis \( \{L_{p}\} \) of \( \mathfrak{sl}(n+1) / \mathfrak{so}(n+1) \). The corresponding functions \( \xi_{p}(z) = (L_{p})_{ab}z^{a}z^{b} \) on the sphere form a basis of quadratic spherical harmonics. Suppose \( \Lambda(t) = \sum_{p}c^{p}(t)L_{p} = [c]^{T}[L] \), where square brackets denote matrices and vectors. From \cref{isotropicZ}, the necessary and sufficient condition for the flow to stay at the critical point is \( D_{L_{q}}\tilde{Z}_{1} = 0 \) for each \( q \in \{1,\dots, n(n+3) / 2\} \), or
\[ \fint_{\Sp^{n}}Q_{n}[\sigma]^{1-\alpha}\xi_{q} = 0 \quad \forall q. \]
Differentiating with respect to time, we get \begin{align*}
    (1-\alpha)\fint_{\Sp^{n}}Q_{n}^{1-\alpha}(A^{-1})^{ij}A[-Q_{n}^{-\alpha} + c^{p}(\sigma\xi_{p} + \frac{1}{2}\left< \nabla \sigma, \nabla \xi_{p} \right>)]_{ij}\xi_{q} & = 0 \quad \forall q
\end{align*}
Note that, with \( \gamma = 1-\alpha(n+2) \),
\[ D_{L_{p}}D_{L_{q}}\tilde{Z}_{1} = -\gamma(1-\alpha)\fint_{\Sp^{n}}Q_{n}^{1-\alpha}(A^{-1})^{ij}A\Big[\sigma \xi_{p} + \frac{1}{2}\left< \nabla \sigma, \nabla \xi_{p} \right>\Big]_{ij}\xi_{q}. \]
By \cref{DAQ}, \( Q_{n}(A^{-1})^{ij}A[\sigma\xi_{p} + \frac{1}{2}\langle \nabla\sigma,\nabla\xi_{p}\rangle]_{ij} = -(n+2)Q_{n}\xi_{p} + \frac{1}{2}\langle \nabla Q_{n},\nabla\xi_{p}\rangle \). Integrating the gradient term by parts and using \( \Delta\xi_{p} = -2(n+1)\xi_{p} \), we get
\begin{equation}
    D_{L_{p}}D_{L_{q}}\tilde{Z}_{1} = \gamma\fint_{\Sp^{n}}Q_{n}^{1-\alpha}\Big( \gamma\xi_{p}\xi_{q} + \frac{1}{2}\langle \nabla\xi_{p},\nabla\xi_{q}\rangle \Big). \label{HessianZpq}
\end{equation}
In particular, the matrix \( D^{2}\tilde{Z}_{1} \) is symmetric, and by \cref{isotropicZ} it is positive-definite as well. Consider the vector \( M = [M^{q}] \) with entries  
\[  M^{q} = (1-\alpha)\gamma\fint_{\Sp^{n}}Q_{n}^{1-\alpha}(A^{-1})^{ij}A[-Q_{n}^{-\alpha}]_{ij}\xi_{q} \]
Now the constraint equation is \begin{equation}
    [M] = [D^{2}\tilde{Z}][c] \label{constrainteq}
\end{equation}
and the coefficients \( c^{p}(t) \) are given by
\[ c^{p}(t) = [D^{2}\tilde{Z}]^{-1}[M]. \] 
For a fixed \( t \), we can choose the basis \( \{L_{p}\} \) so that \( D^{2}\tilde{Z}_{1} \) is diagonal, and then by \cref{HessianZpq}
\begin{equation}
    c^{p}(t) = \frac{M^{p}}{D^{2}\tilde{Z}_{1}(L_{p},L_{p})} = \frac{(1-\alpha)\fint_{\Sp^{n}}Q_{n}^{1-\alpha}(A^{-1})^{ij}A[-Q_{n}^{-\alpha}]_{ij}\,\xi_{p}}{\fint_{\Sp^{n}}Q_{n}^{1-\alpha}\left( \gamma \xi_{p}^{2} +  \frac{1}{2}|\nabla \xi_{p}|^{2}\right)}.\label{cp}
\end{equation}

The short-time existence and uniqueness of ABGCF can be established from standard parabolic PDE theory.
\begin{theorem}\label{shorttime}
    Let \( \mathcal{M}_{0 }\) be a smooth, closed, locally uniformly convex hypersurface in \( \R^{n+1} \) in \( Z_{1} \)-isotropic position. Then there exists a unique, smooth solution of \cref{flow1} for \( t \in [0,T) \) for some \( T > 0 \).
\end{theorem}
\begin{proof}
    With respect to a fixed basis \( \{L_{p}\} \), the global term \( \Lambda(t) \) lives in the finite-dimensional space \( \mathfrak{sl}(n+1) / \mathfrak{so}(n+1) \) and is determined by the constraint equation \cref{constrainteq}. Notice that by integration by parts, the quantity \( M \) is 3rd-order in \( \sigma \) and the Hessian \( D^{2}\tilde{Z}_{1} \) is 2nd-order in \( \sigma \). We write the ABGCF operator as 
    \[ \mathcal{F}(\sigma) = - Q_{n}[\sigma]^{-\alpha} + c^{p}(\sigma)\left( \sigma \xi_{p} + \frac{1}{2}\left< \nabla \sigma, \nabla \xi_{p} \right> \right), \]
    with linearization 
    \begin{align*}
        D \mathcal{F}(\sigma_{0})[h] & = \alpha Q_{n}[\sigma_{0}]^{-\alpha}(A^{-1})^{ij} \left( \nabla_{i}\nabla_{j}h + g_{ij}h \right) + Dc^{p}(\sigma_{0})[h]\left( \sigma_{0} \xi_{p} + \frac{1}{2}\left< \nabla \sigma_{0}, \nabla \xi_{p} \right> \right) \\
        & \quad +  c^{p}(\sigma_{0})\left( \xi_{p}h + \frac{1}{2}\left< \nabla h, \nabla \xi_{p} \right> \right).
    \end{align*}
    Let \( h^{k+\beta}(\Sp^{n}) \) denote the little H\"older space, which is the completion of \( C^{\infty}(\Sp^{n}) \) in the H\"older space \( C^{k,\beta}(\Sp^{n}) \). Let \( X = h^{1+\beta}(\Sp^{n}) \), \( D = h^{3+\beta}(\Sp^{n}) \) with \( \mathcal{O} = \{\sigma \in D : A[\sigma] > 0\} \subset D\) an open set. Then \( \mathcal{F} \in  C^{1}(\mathcal{O},X) \) and the non-local terms form a bounded map from \( h^{3+ \frac{\beta}{2}}(\Sp^{n}) \to h^{1 + \beta}(\Sp^{n}) \) so by \cite[Theorem 2.4.1]{lunardi2012analytic} the total linearization is a sectorial operator. The short-time existence and uniqueness of the flow then follows from \cite[Theorem 8.1.1]{lunardi2012analytic}.

    Finally, the coefficients \( c^{p} \) are chosen so that the time derivative of \( \int Q_{n}^{1-\alpha}\xi_{q}d \mu \) vanishes for every \( q \). Since these integrals vanish at \( t=0 \), the solution stays in \( Z_{1} \)-isotropic position.
\end{proof}
Because the hypersurface stays at the affine critical point of the entropy, the affine corrections don't change the monotonicity of the entropy. Although the proof easily follows from the reference, we will need the exact expression of the time derivative of entropy later, so we provide a short proof of monotonicity here.
 \begin{proposition}\label{Monotonicity}
    Let \( X :M^{n} \times [0,T) \to \R^{n+1}\) be a smooth solution of ABGCF with locally uniformly convex initial hypersurface. Then the entropy \( Z_{1}(\mathcal{M}_{t}) \) is non-decreasing function of time, and strictly increasing unless the solution is a self-similar shrinking solution.
\end{proposition}
\begin{proof}
    We can write \( \log Z_{1} \) as 
    \[ \log Z_{1} = n(\alpha-1)\log \fint_{\Sp^{n}} \sigma Q_{n} d \mu + (n+1) \log \fint_{\Sp^{n}}Q_{n}^{1-\alpha} d \mu. \]
    It follows that under ABGCF, the time derivative is \begin{align*}
        \frac{d}{dt}\log Z_{1} & = \frac{n(\alpha-1)}{\fint_{\Sp^{n}} \sigma Q_{n}} \fint_{\Sp^{n}}(n+1)\partial_{t}\sigma Q_{n}+ \frac{(n+1)(1-\alpha)}{\fint_{\Sp^{n}}Q^{1-\alpha}}\fint_{\Sp^{n}}Q_{n}^{-\alpha}\partial_{t}Q_{n} \\
        & = \frac{n(\alpha-1)(n+1)}{\fint_{\Sp^{n}} \sigma Q_{n}} \fint_{\Sp^{n}} -Q_{n}^{1-\alpha} + \frac{n(1-\alpha)(n+1)}{\fint_{\Sp^{n}}Q_{n}^{1-\alpha}} \fint_{\Sp^{n}}Q_{n}^{-\alpha}Q[-Q_{n}^{-\alpha}, \sigma, \dots, \sigma],
    \end{align*}
    where we are using the fact that volume isn't changed by affine perturbations and the entropy also doesn't change because the hypersurface is at the affine critical point. The positivity now follows from Brunn-Minkowski inequality. Moreover, we record the derivative here for later use
    \begin{align}
        \frac{d}{dt}Z_{1} & = \frac{n(1-\alpha)(n+1)Z_{1}}{\fint \sigma Q_{n}\fint Q_{n}^{1-\alpha}}\left[\left(\fint_{\Sp^{n}}Q_{n}^{1-\alpha}\right)^{2} - \fint \sigma Q_{n} \fint Q_{n}^{-\alpha}Q[Q_{n}^{-\alpha}, \sigma, \dots, \sigma]\right] \\
        & = \frac{n(1-\alpha)(n+1)Z_{1}}{\fint \sigma Q_{n}\fint Q_{n}^{1-\alpha}} \left[\left(\fint_{\Sp^{n}}Q_{n}^{1-\alpha}\right)^{2} + \fint \sigma Q_{n} \fint \nabla_{i}Q_{n}^{-\alpha}\nabla_{j}Q_{n}^{-\alpha}\dot{Q}^{ij} - \fint \sigma Q_{n}\fint Q^{-2 \alpha}Q[1,\sigma,\dots, \sigma]\right]\label{dtZ1}
    \end{align}
    which is non-negative by Brunn-Minkowski inequality for \( f = Q_{n}^{-\alpha} \) with equality holding if and only if \( Q_{n}^{-\alpha} = c \sigma + \left< p,z \right> \) for some \( c \in \R \) and \( p \in \R^{n+1} \) which describes self-similar solutions.
\end{proof}
\begin{proposition}[Evolution equations] \label{evolution}
    Let \( \sigma : \Sp^{n} \times [0,T) \to \R \) be a smooth solution of \cref{flow2}. The following evolution equations hold:
    \begin{equation}
    \partial_{t}Q_{n}[\sigma]^{-\alpha} = \mathcal{L}Q_{n}^{-\alpha} + \alpha H Q_{n}^{-2 \alpha} + \alpha(n+2)Q_{n}^{-\alpha}\xi + \frac{1}{2}\left< \nabla Q_{n}^{-\alpha}, \nabla \xi \right>, \label{evolQ}
    \end{equation}
    and
    \begin{align*}
    \partial_{t} A[\sigma]_{ij} & = \mathcal{L}A_{ij} + (n \alpha - 1)(Q_n)^{-\alpha}g_{ij}  - \alpha H (Q_n)^{-\alpha}A_{ij} \\
    & \quad - \alpha^{2}(Q_n)^{-2-\alpha}\nabla_{i}Q_n\nabla_{j}Q_n - \alpha Q_n^{-\alpha}(A^{-1})^{kp}(A^{-1})^{lq}\nabla_{i}A_{kl}\nabla_{j}A_{pq} \\
    & \quad + A_{ij}\xi + \frac{1}{2}\left< \nabla A_{ij}, \nabla \xi \right> + \frac{1}{2}g^{pq}A_{ip}\nabla_{q}\nabla_{j}\xi + \frac{1}{2}g^{pq}A_{jp}\nabla_{q}\nabla_{i}\xi,
\end{align*}
where \( \mathcal{L} = \alpha Q_n^{-\alpha}(A^{-1})^{ij}\nabla_{i}\nabla_{j} \) is an elliptic operator and \( H = \operatorname{tr}A[\sigma]^{-1} \) is the mean curvature. 
\end{proposition}
\begin{proof}
    For evolution of speed, we add the contribution from affine term \( \xi \) using \cref{DAQ}  to get \begin{align*}
        \partial_{t}Q_{n}^{-\alpha} & = \partial_{t}\det A[\sigma]^{-\alpha} \\
        & = -\alpha Q_{n}^{-\alpha}(A^{-1})^{ij}A[-Q_{n}^{-\alpha} + \sigma \xi + \frac{1}{2}\left< \nabla \sigma, \nabla \xi \right>]_{ij} \\
        & = \alpha Q_{n}^{-\alpha}(A^{-1})^{ij}\left( \nabla_{i}\nabla_{j}Q_{n}^{-\alpha} + g_{ij}Q_{n}^{-\alpha} \right) + \alpha(n+2)Q_{n}^{-\alpha}\xi + \frac{1}{2}\left< \nabla Q_{n}^{-\alpha}, \nabla \xi \right>.
    \end{align*}
    The evolution of individual components \( A_{ij} \) under \( \alpha \)-Gauss curvature flow can be found in \cite[Section 6]{andrews2022extrinsic} which make up the terms in the first two lines. We calculate the contribution from the global term \( \xi \) which is the last line.  Note that  \( z \) is the position vector as well as the normal on \( \Sp^{n} \) and the second fundamental form is equal to the metric  so \( \nabla_{i}\nabla_{j}z = - h_{ij}z = - g_{ij}z \). It follows that \[ \nabla_{i}\nabla_{j}z^{a} =  \nabla_{i}\nabla_{j} \langle e^{a},z \rangle= \langle e^{a}, \nabla_{i}\nabla_{j}z \rangle = -g_{ij}\langle e^{a},z \rangle = -g_{ij}z^{a}. \]
    For \( \xi = L_{ab}z^{a}z^{b} \) with  \( L \) symmetric and \( \operatorname{tr}(L) = 0 \), we calculate   
    \begin{align*}
        \nabla_{j}\nabla_{q}\xi & = L_{ab} \nabla_{j}\nabla_{q}(z^{a}z^{b})   \\
        & = L_{ab}\left( z^{b}\nabla_{j}\nabla_{q}z^{a} + z^{a}\nabla_{j}\nabla_{q}z^{b} + \nabla_{j}z^{a}\nabla_{q}z^{b} + \nabla_{q}z^{a}\nabla_{j}z^{b} \right) \\
        & = L_{ab}\left( -2g_{jq}z^{a}z^{b} + 2\nabla_{j}z^{a}\nabla_{q}z^{b} \right).
    \end{align*}
     Thus, \begin{align*}
        \nabla_{i}\nabla_{j}\nabla_{q}\xi & = - 2g_{jq}\nabla_{i}\xi + 2L_{ab}\nabla_{i}\nabla_{j}z^{a}\nabla_{q}z^{b} + 2L_{ab}\nabla_{j}z^{a}\nabla_{i}\nabla_{q}z^{b} \\
        & = - 2 g_{jq}\nabla_{i}\xi - g_{ij}\nabla_{q}\xi - g_{iq}\nabla_{j}\xi.
        \end{align*}
        Moreover, by Ricci commutator identity on \( \Sp^{n} \), \begin{align*}
        \nabla_{i}\nabla_{j}\nabla_{p}\sigma  = \nabla_{i}\nabla_{p}\nabla_{j}\sigma & = \nabla_{p}\nabla_{i}\nabla_{j}\sigma + R_{ipj}^{r}\nabla_{r}\sigma \\
        & = \nabla_{p}\nabla_{i}\nabla_{j}\sigma + (g_{ij}\delta^{r}_{p} - g_{pj}\delta^{r}_{i})\nabla_{r}\sigma \\
        & = \nabla_{p}\nabla_{i}\nabla_{j}\sigma + g_{ij}\nabla_{p}\sigma - g_{pj}\nabla_{i}\sigma.
        \end{align*}
        The Hessian of the inner product can now be calculated as 
        \begin{align*}
        \nabla_{i}\nabla_{j}\left< \nabla \sigma, \nabla \xi \right>  & = g^{pq}\nabla_{i}\nabla_{j}\nabla_{p}\sigma \nabla_{q}\xi + g^{pq}\nabla_{j}\nabla_{p}\sigma \nabla_{i}\nabla_{q}\xi + g^{pq}\nabla_{i}\nabla_{p}\sigma \nabla_{j}\nabla_{q}\xi + g^{pq}\nabla_{p}\sigma \nabla_{i}\nabla_{j}\nabla_{q}\xi \\
        & = g^{pq}\left( \nabla_{p}\nabla_{i}\nabla_{j}\sigma + g_{ij}\nabla_{p}\sigma - g_{pj}\nabla_{i}\sigma \right)\nabla_{q}\xi +  g^{pq}\left( A[\sigma]_{jp} - g_{jp}\sigma \right)\nabla_{i}\nabla_{q}\xi \\
        & \quad + g^{pq}\left( A[\sigma]_{ip}-g_{ip}\sigma \right)\nabla_{j}\nabla_{q}\xi - g^{pq}\nabla_{p}\sigma \left( g_{ij}\nabla_{q}\xi + g_{iq}\nabla_{j}\xi + 2 g_{jq}\nabla_{i}\xi \right) \\
        & = g^{pq}\nabla_{p}(A[\sigma]_{ij}-g_{ij}\sigma)\nabla_{q}\xi + g_{ij}\left< \nabla \sigma, \nabla \xi \right> - \nabla_{i}\sigma \nabla_{j}\xi + g^{pq}A[\sigma]_{jp}\nabla_{i}\nabla_{q}\xi - \sigma \nabla_{i}\nabla_{j}\xi \\
        & \quad + g^{pq}A[\sigma]_{ip}\nabla_{j}\nabla_{q}\xi - \sigma\nabla_{j}\nabla_{i}\xi - g_{ij}\left< \nabla \sigma, \nabla \xi \right> - \nabla_{i}\sigma \nabla_{j}\xi - 2 \nabla_{j}\sigma \nabla_{i}\xi \\
        & = \left< \nabla A[\sigma]_{ij},\nabla \xi \right> - g_{ij}\left< \nabla \sigma, \nabla \xi \right> + g^{pq}A[\sigma]_{jp}\nabla_{i}\nabla_{q}\xi + g^{pq}A[\sigma]_{ip}\nabla_{j}\nabla_{q}\xi - 2\sigma \nabla_{i}\nabla_{j}\xi \\
        & \quad - 2\nabla_{i}\sigma \nabla_{j}\xi - 2\nabla_{j}\sigma \nabla_{i}\xi.
        \end{align*}
        The contribution from affine term using \cref{DAsigma} is \begin{align*}
        D_{\xi}A[\sigma]_{ij} & = A[\sigma \xi + \frac{1}{2}\left< \nabla \sigma, \nabla \xi \right>]_{ij} \\
        & = A[\sigma]_{ij}\xi + \nabla_{i}\sigma \nabla_{j}\xi + \nabla_{j}\sigma \nabla_{i}\xi + \sigma \nabla_{i}\nabla_{j}\xi + \frac{1}{2} \biggl(\left< \nabla A[\sigma]_{ij},\nabla \xi \right> - g_{ij}\left< \nabla \sigma, \nabla \xi \right> + g^{pq}A[\sigma]_{jp}\nabla_{i}\nabla_{q}\xi \\
        & \quad + g^{pq}A[\sigma]_{ip}\nabla_{j}\nabla_{q}\xi - 2\sigma \nabla_{i}\nabla_{j}\xi - 2\nabla_{i}\sigma \nabla_{j}\xi - 2\nabla_{j}\sigma \nabla_{i}\xi \biggr) + \frac{1}{2}g_{ij}\left< \nabla \sigma, \nabla \xi \right> \\
        & = A[\sigma]_{ij}\xi + \frac{1}{2}\left< \nabla A[\sigma]_{ij}, \nabla \xi \right> + \frac{1}{2}g^{pq}A[\sigma]_{jp}\nabla_{i}\nabla_{q}\xi + \frac{1}{2}g^{pq} A[\sigma]_{ip} \nabla_{j}\nabla_{q}\xi.
    \end{align*}
\end{proof}
The flow has the following symmetries. In particular, \( Z_{1} \)-entropy is translation invariant, but a translate of a solution remains a solution only if the translation vector moves with the affine term.
\begin{lemma}[Symmetries of ABGCF]\label{translationlemma}
    Let \( X \) be a solution of ABGCF with coefficient \( \Lambda(t) \).
    \begin{enumerate}
        \item For \( O \in O(n+1) \), \( OX \) is a solution of ABGCF with coefficient \( O\Lambda O^{T} \).
        \item For \( r > 0 \), the rescaling \( X_{r}(\cdot,t) = rX\big( \cdot, r^{-(1+n\alpha)}t \big) \) is a solution of ABGCF with coefficient
        \[ \Lambda_{r}(t) = r^{-(1+n\alpha)}\Lambda\big( r^{-(1+n\alpha)}t \big). \]
        \item For a \( C^{1} \) curve \( q \) in \( \R^{n+1} \), \( X - q \) is a solution of ABGCF with coefficient \( \Lambda \) if and only if \( \frac{dq}{dt} = \Lambda q \).
    \end{enumerate}
\end{lemma}
\begin{proof}
    \( Z_{1} \)-isotropy is invariant under orthogonal transformations, dilations and translations, and the coefficient of ABGCF is uniquely determined by the solution, since \( D^{2}\tilde{Z}_{1} \) is invertible. 
    
    (1) We have \( \partial_{t}(OX) = -K^{\alpha}O\nu + O\Lambda O^{T}(OX) \), and \( O\nu \) is the normal of \( OX \). 
    
    (2) The Gauss curvature of \( X_{r} \) is \( r^{-n}K \), so \( \partial_{t}X_{r} = r^{-n\alpha}\big( -K^{\alpha}\nu + \Lambda X \big) = -K_{r}^{\alpha}\nu + \Lambda_{r}X_{r} \). 
    
    (3) We have \( \partial_{t}(X - q) = -K^{\alpha}\nu + \Lambda(X - q) + \Lambda q - q' \).
\end{proof}
\section{Geometry at \( Z_{1} \)-isotropic position} \label{geometryatZ1}

In this section we establish entropy, isoperimetric and diameter estimates for a hypersurface at the \( Z_{1} \)-isotropic position. Throughout the rest of the paper \( C \) will denote a constant depending on the parameters set in the hypothesis which may change in different proofs.
\begin{proposition}\label{Z1upperbound}
    Let \( \mathcal{M} = \partial \Omega \) be a smooth hypersurface at \( Z_{1} \)-isotropic position. The Gaussian entropy \( Z_{1} \) satisfies an upper bound, 
    \[ Z_{1}(\Omega) \le C_{\text{iso}}(n).  \]
\end{proposition}
\begin{proof}
    Let \( J \) be the John ellipsoid of \( \mathcal{M} \), so there exists a \( q \in \Omega \) such that \( J \subset \Omega - q \subset (n+1)J  \). Since entropy is scale invariant we can assume that \( V(\mathcal{M}) = V(B^{n+1}) \). Let \( T \in \operatorname{SL}(n+1) \) such that \( TJ = B_{r}(0) \) is a ball of radius \( r \) centered at the origin. By inclusion property of John ellipsoid we have 
\[ B_{r}(0) \subset T(\Omega - q)\subset B_{(n+1)r}(0). \]
The volume constraint now implies \( r \le 1 \le (n+1)r \) which also gives upper and lower bound on the support function so \( \frac{1}{n+1} \le r \le \sigma_{T(\Omega-q)}\le (n+1)r \le (n+1) \). Now, 
\[ Z_{\frac{1}{\alpha}}(\sigma_{T(\Omega-q)}) = \left( \fint_{\Sp^{n}}(\sigma_{T(\Omega-q)})^{1-\frac{1}{\alpha}} \right)^{n+1} \le (n+1)^{(\frac{1}{\alpha}-1)(n+1)}. \]
At isotropic position, by H\"older inequality 
\[ Z_{1}(\sigma_{\Omega}) = Z_{1}(\sigma_{\Omega-q}) \le Z_{1}(\sigma_{T(\Omega-q)}) \le Z_{\frac{1}{\alpha}}(\sigma_{T(\Omega-q)})^{\alpha} \le (n+1)^{(n+1)(1-\alpha)} \le (n+1)^{n+1}, \]
which completes the proof.
\end{proof}

\begin{proposition}\label{isoperimetric}
    Let \( \mathcal{M} \) be a hypersurface in \( Z_{1} \)-isotropic position. Then there exists a positive constant \( C(n,\alpha) \) such that 
    \[ \frac{\rho_{+}(\mathcal{M})}{\rho_{-}(\mathcal{M})} \le CZ_{1}(\mathcal{M})^{-\frac{1}{1-\alpha}}, \]
    where \( \rho_{+} \) and \( \rho_{-} \) denote the circumradius and inradius respectively.
\end{proposition}
\begin{proof}
    We use projected area to bound the isoperimetric ratio like Theorem 4 in \cite{andrews2000motion}. Let \( z_{0} \in \Sp^{n} \) be such that the projected area of \( \mathcal{M} \) onto \( z_{0}^{\perp} \) plane is minimal. The projected area of \( \mathcal{M} \) onto \( z_{0}^{\perp} \) plane is given by 
    \[ 2A_{z_{0}} = \int_{\Sp^{n}}|\langle z,z_{0} \rangle|Q_{n}d \mu. \]
    From H\"older inequality and isotropic condition, we get 
  \begin{align*}
       2A_{z_{0}} & \ge \int_{\Sp^{n}}| \langle z,z_{0} \rangle|^{\frac{2}{1-\alpha}}Q_{n}d \mu \\
       & \ge \left( \int_{\Sp^{n}} |\langle z,z_{0}\rangle|^{2}Q_{n}^{1-\alpha}d \mu\right)^{\frac{1}{1-\alpha}}|\Sp^{n}|^{ \frac{-\alpha}{1-\alpha}} \\
       & = \left( \frac{|\Sp^{n}|}{n+1} \fint_{\Sp^{n}}Q_{n}^{1-\alpha}d \mu \right)^{\frac{1}{1-\alpha}}|\Sp^{n}|^{ \frac{-\alpha}{1-\alpha}} = |\Sp^{n}|(n+1)^{\frac{-1}{1-\alpha}}V^{ \frac{n}{n+1}}Z_{1}^{\frac{1}{(1-\alpha)(n+1)}}.
  \end{align*}
  Using geometric bounds from Theorem 4 in \cite{andrews2000motion} (or Corollary 17.8 in \cite{andrews2022extrinsic}), there exists \( C(n) \) such that 
  \[ \frac{\rho_{+}}{\rho_{-}} \le C(n) \frac{V^{n}}{A_{z_{0}}^{n+1}}  \] which along with the lower bound on \( A_{z_{0}} \) gives the desired result.
\end{proof}
Since entropy is non-decreasing, the above proposition bounds the isoperimetric ratio throughout the flow, preventing the hypersurface from degenerating in one direction. 
\begin{corollary}
    Let \( \{\mathcal{M}_{t}\} \) be a solution of ABGCF with \( \mathcal{M}_{0} \) at \( Z_{1} \)-isotropic position. Then there exists a constant \( C=C(n,\alpha,Z_{1}(\mathcal{M}_{0})) \) such that 
    \[ \frac{\rho_{+}(\mathcal{M}_{t})}{\rho_{-}(\mathcal{M}_{t})} \le C. \]\label{isoperimetricbounded}
\end{corollary}
    With the above corollary, the control of volume implies the control of inradius, circumradius and diameter for ABGCF.
\begin{proposition}
    Let \( \{\mathcal{M}_{t}\} \) be a solution of ABGCF with \( \mathcal{M}_{0} \) at \( Z_{1} \)-isotropic position. Then there exists a constant \( C = C(n,\alpha,Z_{1}(\mathcal{M}_{0})) \) such that 
    \[V(\mathcal{M}_{0})^{ \frac{1+n \alpha}{1+n}}- Ct \le V(\mathcal{M}_{t})^{ \frac{1+n \alpha}{1+n}} \le V(\mathcal{M}_{0})^{ \frac{1+n \alpha}{1+n}} - C^{-1}t \]
    and \[ C^{-1}\left( V(\mathcal{M}_{0})^{ \frac{1+n \alpha}{1+n}}- Ct  \right)^{\frac{1}{1+n \alpha}} \le \rho_{-}(\mathcal{M}_{t}) \le \rho_{+}(\mathcal{M}_{t}) \le C\left( V(\mathcal{M}_{0})^{ \frac{1+n \alpha}{1+n}}- C^{-1}t  \right)^{\frac{1}{1+n \alpha}}. \]\label{Volumeandgeometry}
\end{proposition}
\begin{proof}
    The derivative of volume under the flow isn't affected by the affine term, so we get
    \begin{align*}
        \frac{d}{dt} \frac{V(\mathcal{M}_{t})}{V(B^{n+1})} & = - (n+1)\fint_{\Sp^{n}}Q_{n}^{1-\alpha}\,d\mu \\
        & = - (n+1)Z_{1}(\mathcal{M}_{t})^{\frac{1}{n+1}}\left( \frac{V(\mathcal{M}_{t})}{V(B^{n+1})} \right)^{ \frac{n(1-\alpha)}{1+n}}.
    \end{align*}
    The entropy is bounded above by a dimensional constant and below by \( Z_{1}(\mathcal{M}_{0}) \), so we get the desired inequality. By the isoperimetric ratio bound we also know that \( \rho_{+}(\mathcal{M}_{t}) \sim \rho_{-}(\mathcal{M}_{t}) \sim V(\mathcal{M}_{t})^{\frac{1}{1+n}} \) which finishes the proof.
\end{proof}
\section{Upper speed bound}\label{sec:upperspeed}

In this section we prove an \textit{a priori} upper bound on speed as long as the inradius is bounded.  It will be essential to estimate the global coefficient \( \xi(t) \) of \cref{flow2} in order to get curvature estimates. We write \( |\xi(t)| = \sup_{z \in \Sp^{n}}|\xi(z,t)| \), which is the operator norm of \( \Lambda(t) \). Since \( \nabla \xi = 2(\Lambda z)^{T} \) and \( \nabla_{i}\nabla_{j}\xi = 2\langle \Lambda e_{i},e_{j} \rangle - 2\xi g_{ij} \), we have \( |\nabla \xi| + |\nabla^{2}\xi| \le C(n)|\xi(t)| \).

\begin{lemma}
    Let \( \sigma : \Sp^{n} \times [0,T_{0}] \to \R \) be a hypersurface flowing by ABGCF with inradius and circumradius bounds \( r_{0} \le \rho_{-}(\sigma_{t}) \le \rho_{+}(\sigma_{t}) \le R_{0} \) for \( t \in [0,T_{0}] \), then there exists a constant \( C(n,\alpha) > 0\) such that 
    \[ |\xi(t)|\le C(n,\alpha)\left( R_{0}^{n-1}r_{0}^{n(\alpha-1)}Z_{1}(\sigma_{0})^{-\frac{1}{n+1}} \sup_{z \in \Sp^{n}}Q_{n}(z,t)^{-2 \alpha} + \frac{d}{dt}\log Z_{1}(\sigma(t))\right). \]\label{xilemma}
\end{lemma}
\begin{proof}
    Since the scale of the hypersurface remains comparable during the flow, and we can bound volume, surface area and quermassintegrals above and below by constants depending on \( r_{0},R_{0} \) and \( n \).

    Fix a basis \( \{L_{p}\} \) of \( \mathfrak{sl}(n+1) / \mathfrak{so}(n+1) \) such that the hessian \( D^{2}Z_{1} \) with respect to affine perturbation is diagonal in this basis. Using the formula from \cref{isotropicZ} the Hessian of \( \tilde{Z}_{1} \) in the affine direction is
    \begin{align}
        D^{2}\tilde{Z}_{1}(L_{p},L_{p}) = \frac{d^{2}}{ds^{2}}\bigg|_{s=0}\tilde{Z}_{1}(\sigma_{s}) & = \fint_{\Sp^{n}}Q_{n}[\sigma]^{1- \alpha} \left(\gamma^{2}\left< y,L_{p}y \right>^{2} + 2 \gamma \left(|L_{p}y|^{2} - \left< y,L_{p}y \right>^{2}\right)\right)d \mu(y) \nonumber \\
        & = \fint_{\Sp^{n}}Q_{n}[\sigma]^{1- \alpha}\left( \gamma(2- \gamma) \left( \left< y,L_{p}^{2}y \right> - \left< y,L_{p}y \right>^{2} \right) + \gamma^{2} \left< y,L_{p}^{2}y \right> \right)d \mu(y) \nonumber \\
        & \ge \gamma^{2} \fint_{\Sp^{n}}Q_{n}[\sigma]^{1- \alpha}\left< y,L_{p}^{2}y \right>d \mu(y) = \frac{\gamma^{2} \operatorname{tr}(L_{p}^{2})}{n+1}\fint_{\Sp^{n}}Q_{n}[\sigma]^{1-\alpha}, \label{HessianZ}
    \end{align}
    where we used the fact that for \( \beta = 1 \), \( \gamma = 1-\alpha(n+2)<2 \) and the hypersurface is at isotropic \( Z_{1} \)-position for the last equality. From \cref{HessianZ}, the denominator of \( c^{p} \) satisfies the inequality \( \gamma^{-1} D^{2}\tilde{Z}_{1}(L_{p},L_{p}) \ge \gamma \operatorname{tr}(L_{p}^{2})\tilde{Z}_{1}(\sigma) \ge C_{1}(n,\alpha)\fint_{\Sp^{n}}Q_{n}^{1-\alpha} \).
   
    To control the numerator, we use integration by parts (which is valid since \( \dot{Q}^{ij} \) is divergence free) and regroup some terms, \begin{align*}
    \fint_{\Sp^{n}}\xi_p  Q_n^{-\alpha}Q[Q_n^{-\alpha}, \sigma,\dots, \sigma] & = \fint_{\Sp^{n}} \xi_p  Q_n^{-2 \alpha}Q[1,\sigma,\dots, \sigma] - \fint_{\Sp^{n}} \nabla_{i}(\xi_p  Q_n^{-\alpha})\nabla_{j}Q_n^{-\alpha}\dot{Q}^{ij} \\
    & = \fint_{\Sp^{n}} \xi_p  Q_n^{-2 \alpha}Q[1,\sigma,\dots, \sigma] - \fint_{\Sp^{n}} \xi_p  \nabla_{i}Q_n^{-\alpha}\nabla_{j}Q_n^{-\alpha}\dot{Q}^{ij} - \frac{1}{2}\fint_{\Sp^{n}} \nabla_{i}\xi_p  \nabla_{j}Q_n^{-2 \alpha}\dot{Q}^{ij} \\
    & = \fint_{\Sp^{n}} \xi_p  Q_n^{-2 \alpha}Q[1,\sigma,\dots, \sigma] - \fint_{\Sp^{n}} \xi_p  \nabla_{i}Q_n^{-\alpha}\nabla_{j}Q_n^{-\alpha}\dot{Q}^{ij} \\
    & \quad + \frac{1}{2}\fint_{\Sp^{n}} Q_n^{-2 \alpha}Q[\xi_p , \sigma, \dots, \sigma] - \frac{1}{2} \fint_{\Sp^{n}}\xi_p  Q_n^{-2 \alpha}Q[1, \sigma,\dots, \sigma] \\
    & = \frac{1}{2}\fint_{\Sp^{n}} \xi_p  Q_n^{-2 \alpha}Q[1,\sigma,\dots, \sigma] - \fint_{\Sp^{n}} \xi_p  \nabla_{i}Q_n^{-\alpha}\nabla_{j}Q_n^{-\alpha}\dot{Q}^{ij} \\
    & \quad + \frac{1}{2}\fint_{\Sp^{n}} Q_n^{-2 \alpha}Q[\xi_p , \sigma, \dots, \sigma].
\end{align*}
There exists a constant \( C_{2}(n) \) such that \( |\xi_{p}| \le C_{2}(n) \) and \( -C_{2}(n)g_{ij} \le A[\xi_{p}]_{ij} \le C_{2}(n)g_{ij} \). It follows that \( |Q[\xi_{p},\sigma,\dots, \sigma]| \le C_{2}(n) Q[1,\sigma,\dots, \sigma] \) since the mixed discriminant satisfies inequality \( Q[f,\sigma,\dots, \sigma] \ge 0 \) for any \( f \) with \( A[f]_{ij} \ge 0 \). Taking modulus and applying triangle inequality, we get \begin{align*}
    \bigg|\fint_{\Sp^{n}}\xi_p  Q_n^{-\alpha}Q[Q_n^{-\alpha}, \sigma,\dots, \sigma]\bigg| & \le C_{2}(n)  \fint_{\Sp^{n}} Q_{n}^{-2 \alpha}Q[1,\sigma,\dots, \sigma] + C_{2}(n) \fint_{\Sp^{n}} \nabla_{i}Q_{n}^{-\alpha}\nabla_{j}Q_{n}^{-\alpha}\dot{Q}^{ij} \\
    & \le 2C_{2}(n)  \fint_{\Sp^{n}} Q_{n}^{-2 \alpha}Q[1,\sigma,\dots, \sigma] \\
    & \quad + C_{2}(n) \left[ \frac{\fint_{\Sp^{n}}Q_{n}^{1-\alpha}}{n(n+1)(1-\alpha)Z_{1}}\frac{d}{dt}Z_{1}(\sigma(t)) - \frac{\left( \fint_{\Sp^{n}}Q_{n}^{1-\alpha} \right)^{2}}{\fint_{\Sp^{n}}\sigma Q_{n}} \right] \\
    & \le 2C_{2}(n)\sup_{z \in \Sp^{n}}Q_{n}(z,t)^{-2 \alpha}\fint_{\Sp^{n}}Q[1,\sigma,\dots, \sigma] + C_{3}(n,\alpha) \frac{\fint_{\Sp^{n}}Q_{n}^{1-\alpha}}{Z_{1}}\frac{d}{dt}Z_{1}
\end{align*}
where we used \cref{dtZ1} to the time derivative of entropy term. Combining the estimates for numerator and denominator gives \begin{align*}
    |\xi(t)| & \le 2 \frac{C_{2}(n)}{C_{1}(n,\alpha)}\sup_{z \in \Sp^{n}}Q_{n}(z,t)^{-2 \alpha} \frac{\fint_{\Sp^{n}}Q[1,\sigma,\dots, \sigma]}{\fint_{\Sp^{n}}Q_{n}^{1-\alpha}} + \frac{C_{3}(n,\alpha)}{C_{1}(n,\alpha)}\frac{1}{Z_{1}}\frac{d}{dt}Z_{1} \\
    & \le C(n,\alpha) \frac{R_{0}^{n-1}}{Z_{1}(\sigma_{0})^{\frac{1}{n+1}}r_{0}^{n(1-\alpha)}}\sup_{z \in \Sp^{n}}Q_{n}(z,t)^{-2 \alpha} + C(n,\alpha)\frac{d}{dt} \log Z_{1},
\end{align*}
where the denominator in the first term is bounded below by \( Z_{1}(\sigma_{0})^{\frac{1}{n+1}}r_{0}^{n(1-\alpha)} \), using the monotonicity of the entropy and \( V(\sigma_{t}) \ge V(B^{n+1})r_{0}^{n+1} \), and the numerator is a quermassintegral bounded above by \( R_{0}^{n-1} \).
\end{proof}

Many of the techniques for curvature estimates in the rest of the paper will be following estimates for \( \alpha \)-GCF with \( \alpha \le 1 /n \) with adjustments by an entropy term to counter the global term. The following theorem provides an upper bound on the speed of the flow.     
\begin{theorem}\label{upperspeed}
    Let \( \sigma : \Sp^{n} \times [0,T_{0}] \to \R \) be a solution of ABGCF with \( r_{0} \le \sigma_{t} \le R_{0} \) for \( t \in [0,T_{0}] \), that is, \( B(0,r_{0}) \subset \Omega_{t} \subset B(0,R_{0}) \). Then there exists a constant \( C(n,\alpha,Z_{1}(\sigma_{0})) > 0 \) such that
    \[  Q_{n}(z,t)^{-\alpha} \le C(n,\alpha,Z_{1}(\sigma_{0})) \frac{R_{0}}{r_{0}}\left( r_{0}^{-n \alpha} + R_{0}^{ \frac{n(n+1)\alpha}{1-n \alpha}}r_{0}^{ \frac{-n^{2}\alpha(1-\alpha)}{1-n \alpha}- \frac{2n \alpha}{1-n \alpha}} + t^{- \frac{n \alpha}{1+n \alpha}} \right) \]
    for every \( z \in \Sp^{n} \) and \( t \in (0,T_{0}] \).
\end{theorem}
\begin{remark}
    By \cref{translationlemma}(3), \cref{upperspeed} applies to \( X - q \) whenever \( q' = \Lambda q \) and \( B(q(t),r_{0}) \subset \Omega_{t} \subset B(q(t),R_{0}) \) for \( t \in [0,T_{0}] \). Since \( Q_{n} \) is translation invariant, the conclusion then holds for \( X \). Such curves \( q \) are constructed in \cref{window}.
\end{remark}
\begin{proof}
    We use Tso's trick \cite{tso1985deforming} to get an upper bound of speed \( Q_{n}^{-\alpha} \). Let \( \varphi = \frac{Q_{n}^{-\alpha}}{\sigma - r_{0} /2} \) which is comparable to \( Q_{n}^{-\alpha} \) since \(  \varphi r_{0} /2 \le Q_{n}^{-\alpha} \le (R_{0}-r_{0} /2)\varphi \le R_{0}\varphi \). The evolution of  \( \varphi \) is given by \begin{align*}
        \partial_{t}\varphi & = \frac{\mathcal{L}Q_{n}^{-\alpha} + \alpha H Q_{n}^{-2 \alpha} + \alpha(n+2)Q_{n}^{-\alpha}\xi + \frac{1}{2}\left< \nabla Q_{n}^{-\alpha}, \nabla \xi \right>}{\sigma - r_{0}/2}  - \frac{Q_{n}^{-\alpha}}{(\sigma - r_{0}/2)^{2}}(-Q_{n}^{-\alpha} + \sigma \xi + \frac{1}{2}\left< \nabla \sigma, \nabla \xi \right>)
    \end{align*}
    where \( \mathcal{L} = \alpha Q_n^{-\alpha}(A^{-1})^{ij}\nabla_{i}\nabla_{j} \). Now \( \mathcal{L}\sigma = n \alpha Q_{n}^{-\alpha} - \alpha H Q_{n}^{-\alpha}\sigma \) and 
    \[ \mathcal{L}\varphi = \frac{\mathcal{L}Q_{n}^{-\alpha}}{\sigma - r_{0} /2} - \frac{Q_{n}^{-\alpha}\mathcal{L}\sigma}{(\sigma-r_{0} /2)^{2}} - 2 \alpha \varphi A^{-1}(\nabla \varphi, \nabla \sigma). \]
    Thus,
    \begin{align*}
        \partial_{t}\varphi & = \mathcal{L}\varphi + \frac{\alpha H Q_{n}^{-2 \alpha}}{\sigma - r_{0}/2} + \frac{Q_{n}^{-\alpha}}{(\sigma - r_{0}/2)^{2}}(n \alpha Q_{n}^{-\alpha} - \alpha H Q_{n}^{-\alpha}\sigma + Q_{n}^{-\alpha}) + 2 \alpha \varphi A^{-1}(\nabla \varphi, \nabla \sigma) \\
        & \quad  + \frac{\alpha(n+2)Q_{n}^{-\alpha}\xi}{\sigma - r_{0}/2} + \frac{1}{2(\sigma - r_{0}/2)}\left< \nabla Q_{n}^{-\alpha}, \nabla \xi \right>  - \frac{Q_{n}^{-\alpha}}{(\sigma - r_{0}/2)^{2}}(\sigma \xi + \frac{1}{2}\left< \nabla \sigma, \nabla \xi \right>) \\
        & = \mathcal{L}\varphi + 2 \alpha \varphi A^{-1}(\nabla \varphi, \nabla \sigma) + (1+ n \alpha - \alpha H r_{0}/2) \varphi^{2} + \left( \alpha(n+2) - 1 - \frac{r_{0}/2}{\sigma - r_{0}/2} \right)\varphi \xi \\
        & \quad + \frac{1}{2}\left< \nabla \varphi, \nabla \xi \right>.
    \end{align*}
    In the presence of global term \( \xi \), we can't directly conclude bounds on the speed. In order to remedy this, we multiply it by large negative power of entropy to control the bad terms and use the estimate from previous lemma. Let \( \psi = Z_{1}^{-\beta}\varphi \), where \( \beta >0 \) will be chosen later. Note that \( \alpha r_{0}H /2  \ge n\alpha r_{0} /2  (Q_{n}^{-\alpha})^{ \frac{1}{n\alpha}} \ge n \alpha r_{0} /2 (Z_{1}^{\beta}\psi r_{0} /2)^{\frac{1}{n\alpha}} \ge C_{1}(n,\alpha) Z_{1}(\sigma_{0})^{ \frac{\beta}{n \alpha}}r_{0}^{1+\frac{1}{n \alpha}} \psi^{\frac{1}{n\alpha}}\). For the linear term, we bound the coefficient as 
    \[ \bigg|\alpha(n+2)-1- \frac{r_{0} / 2}{\sigma - r_{0} /2}\bigg| = 1-\alpha(n+2) + \frac{r_{0} / 2}{\sigma - r_{0} /2}\le 2-\alpha(n+2)= C_{2}(n,\alpha).  \]
    It follows that
    \begin{align*}
        \partial_{t}\psi & = Z_{1}^{-\beta}\partial_{t}\varphi - \beta Z_{1}^{-\beta-1}\varphi \frac{d}{dt}Z_{1} \\
        & \le  \mathcal{L}\psi + 2 \alpha Z_{1}^{\beta}\psi A^{-1}(\nabla \psi, \nabla \sigma) + (1+n \alpha)Z_{1}^{\beta}\psi^{2} + C_{2}|\xi|\psi + \frac{1}{2}\left< \nabla \psi, \nabla \xi \right> \\ 
        & \quad - C_{1}Z_{1}(\sigma_{0})^{\frac{\beta}{n \alpha}+\beta}r_{0}^{1+ \frac{1}{n \alpha}} \psi^{2 + \frac{1}{n\alpha}} - \beta Z_{1}^{-1}\psi \frac{d}{dt}Z_{1} \\
        & \le \mathcal{L}\psi + 2 \alpha Z_{1}^{\beta}\psi A^{-1}(\nabla \psi, \nabla \sigma) + (1+n \alpha)Z_{1}^{\beta}\psi^{2}\\
        & \quad  + C_{2}C\left( R_{0}^{n+1}r_{0}^{n(\alpha-1)}Z_{1}(\sigma_{0})^{-\frac{1}{n+1}} Z_{1}^{2 \beta}\psi_{\max}^{2} + \frac{d}{dt}\log Z_{1}(\sigma(t))\right)\psi \\
        & \quad + \frac{1}{2}\left< \nabla \psi, \nabla \xi \right> - C_{1}Z_{1}(\sigma_{0})^{\frac{\beta}{n \alpha}+\beta}r_{0}^{1+\frac{1}{n \alpha}} \psi^{2 + \frac{1}{n\alpha}} - \beta \psi \frac{d}{dt}\log Z_{1}.
    \end{align*}
    Let \( \beta \) be large enough such that \( \beta(n,\alpha) \ge C_{2}(n,\alpha)C(n,\alpha) \), then we can use \cref{Z1upperbound} to get \( Z_{1}^{\beta} \le C_{3}(n,\alpha) \). Thus,
    \begin{align*}
        \partial_{t}\psi & \le \mathcal{L}\psi + 2 \alpha Z_{1}^{\beta}\psi A^{-1}(\nabla \psi, \nabla \sigma) + C_{3}(1+n \alpha)\psi^{2}+ C_{2}CC_{3}R_{0}^{n+1}r_{0}^{n(\alpha-1)}Z_{1}(\sigma_{0})^{-\frac{1}{n+1}}\psi_{\max}^{3} \\
        & \quad  + \frac{1}{2}\left< \nabla \psi, \nabla \xi \right> - C_{1}Z_{1}(\sigma_{0})^{\frac{\beta}{n \alpha}+\beta}r_{0}^{1+\frac{1}{n \alpha}} \psi^{2 + \frac{1}{n\alpha}}.
    \end{align*}
    At a maximum point of \( \psi \), we have \( \nabla \psi = 0 \) and \( \mathcal{L}\psi \le 0 \), so 
    \begin{align*}
        \frac{d}{dt}\psi_{\max} &\le \left( C_{3}(1+n \alpha) + C_{4}R_{0}^{n+1}r_{0}^{n(\alpha-1)}Z_{1}(\sigma_{0})^{-\frac{1}{n+1}}\psi_{\max} - C_{1}Z_{1}(\sigma_{0})^{\frac{\beta}{n \alpha}+\beta}r_{0}^{1+\frac{1}{n \alpha}} \psi_{\max}^{\frac{1}{n\alpha}}\right)\psi_{\max}^{2} .
    \end{align*}
    Using  \( \frac{1}{n \alpha} > \frac{n+2}{n} > 1 \) we get \( \psi_{\max} \) is bounded above by 
    \begin{align*}
        \psi_{\max}(t) & \le \max \bigg\{C_{5}(n,\alpha)Z_{1}(\sigma_{0})^{-\beta'(n,\alpha)}r_{0}^{-1-n \alpha}, C_{6}(n,\alpha)Z_{1}(\sigma_{0})^{-\beta''(n,\alpha)}R_{0}^{\frac{n(n+1)\alpha}{1-n\alpha}}r_{0}^{ \frac{-n^{2}\alpha(1-\alpha)}{1-n \alpha} - \frac{1+n \alpha}{1-n \alpha}}, \\
        & \quad C_{7}(n,\alpha)Z_{1}(\sigma_{0})^{ \frac{-\beta}{1+n \alpha}}r_{0}^{-1}t^{ \frac{-n \alpha}{1+n \alpha}} \bigg\}.
    \end{align*}
    The speed bound now follows from the relation \( Q_{n}^{-\alpha}(z,t) \le R_{0}Z_{1}^{\beta}(\sigma(t))\psi(z,t) \).
\end{proof}
Absorbing the dependence of \( r_{0} \) and \( R_{0} \) into the constant we can also write 
\[ Q_{n}(z,t)^{- \alpha} \le C(n,\alpha, Z_{1}(\sigma_{0}),r_{0},R_{0}) \left( 1 + t^{ \frac{-n \alpha}{1+n \alpha}} \right). \]
We can now bootstrap the upper speed bound on the global coefficient \( \xi \). The improved estimate is 
\begin{equation}   
    |\Lambda(t)| \le C(n,\alpha,Z_{1}(\sigma_{0}), r_{0},R_{0})\left(1 + t^{ -\frac{2n \alpha}{1+n \alpha}} \right) + C_{2}(n,\alpha,Z_{1}(\sigma_{0}))\frac{d}{dt}\log Z_{1}. \label{xiestimate}
\end{equation}
It follows that \( \xi \) is bounded in \( L^{1}([0,T_{0}]) \) which will be useful in the next section.

\section{Lower speed bound and anisotropic Gauss curvature flow}\label{sec:lowerspeed}
In order to get a lower speed bound, we absorb the affine term using affine transformations to convert it into Gauss curvature flow with time-dependent anisotropy and use spherical barriers to get a lower displacement bound. The displacement bound then gives a lower speed bound by an argument by Smoczyk \cite[Proposition 4]{smoczyk1998starshaped} (also see \cite[Section 9]{andrews2013contracting}).

\begin{lemma} \label{affinepicture}
    Let \( X(t):M^{n} \times [0,T) \to \R^{n+1} \) be a solution of ABGCF. Let \( X^{L}(t) = L(t)X(t) \) where \( L(t) \in \operatorname{GL}(n+1) \) satisfies \[ \frac{d}{dt}L = - L \Lambda, \quad \text{ and } \quad L(0) = \mathrm{Id}. \] Then \( X^{L}(t) \) satisfies 
    \begin{equation}
        \left( \frac{\partial}{\partial t}X^{L} \right)^{\perp} = - \Psi(L(t),\nu^{L})(K^{L})^{\alpha} \label{affineX}
    \end{equation}
    where \( \nu^{L} \) denote the normal of \( X^{L} \) and \( \Psi(L,\nu^{L}) = |L^{T}\nu^{L}|^{1-\alpha(n+2)} \).
\end{lemma}
\begin{proof}
    It is easy to verify that \( L(t) \in \operatorname{SL}(n+1) \) because \( \frac{d}{dt}\det L = \det L \operatorname{tr}(L^{-1}\frac{d}{dt}L) = - \det L \operatorname{tr}(\Lambda) = 0 \). We compute the time derivative of \( X^{L} \):
    \begin{align*}
        \frac{\partial}{\partial t}X^{L} & = \left(\frac{d}{dt}L\right) X + L \frac{\partial}{\partial t}X \\
        & = -L \Lambda X + L(-K^{\alpha}\nu + \Lambda X)  = - K^{\alpha} L \nu.
    \end{align*}
    From \cref{affineQ}, we know that \( \nu^{L} = \frac{L^{-T}\nu}{|L^{-T}\nu|} \) and \( K^{L} = K |L^{-T}\nu|^{-(n+2)} = K|L^{T}\nu^{L}|^{n+2} \). Thus,
    \[ \left( \frac{\partial}{\partial t}X^{L} \right)^{\perp} = - K^{\alpha} \left< L \nu, \nu^{L} \right> = -K^{\alpha}|L^{T}\nu^{L}| = - (K^{L})^{\alpha}|L^{T}\nu^{L}|^{1-\alpha(n+2)}. \]
\end{proof}
Up to a time-dependent reparametrization of \( M^{n} \), the flow \cref{affineX} is a Gauss curvature flow with time-dependent anisotropy. Let \( \sigma^{L} \) be the support function of \( X^{L} \), then \( \sigma^{L} \) satisfies
\begin{equation}
    \partial_{t}\sigma^{L}(z,t) = - \Psi(L(t),z)Q_{n}[\sigma^{L}](z,t)^{-\alpha}\label{supportL}
\end{equation}
with \( \Psi(L,z) = |L^{T}z|^{1-\alpha(n+2)} \). Since \( 1 - \alpha(n+2) \in (0,1) \), \( \|L^{-1}\|^{-1} \le |L^{T}z| \le \|L\| \) gives \( \|L^{-1}\|^{-1} \le \Psi(L,z) \le \|L\| \). Moreover, since \( (L^{-1})' = \Lambda L^{-1} \), we have \( \|L(t)\| + \|L(t)^{-1}\| \le 2\exp\int_{0}^{t}|\Lambda| \).

\begin{corollary}\label{nesting}
    Let \( X \) and \( L \) be as in \cref{affinepicture}, with \( L \) normalised by \( L(t_{1}) = \mathrm{Id} \) at some \( t_{1} \in [0,T) \).
    \begin{enumerate}
        \item The compact convex sets \( L(t)\overline{\Omega}_{t} \) are non-increasing in \( t \in [t_{1},T) \).
        \item The solutions of \( q' = \Lambda q \) on \( [t_{1},T) \) are exactly \( q(t) = L(t)^{-1}q_{1} \), \( q_{1} \in \R^{n+1} \).
        \item There is \( p : [0,T) \to \R^{n+1} \) with \( p' = \Lambda p \) and \( p(t) \in \overline{\Omega}_{t} \) for all \( t \in [0,T) \).
    \end{enumerate}
\end{corollary}
\begin{proof}
    (1) By \cref{affineX}, the hypersurface \( L(t)X(\cdot,t) \) moves with inward normal velocity. 
    
    (2) We have \( \frac{d}{dt}\big( L^{-1}q_{1} \big) = -L^{-1}L'L^{-1}q_{1} = \Lambda L^{-1}q_{1} \), and solutions of the linear equation are unique. 
    
    (3) With \( t_{1} = 0 \), pick \( p^{L} \) in the non-empty intersection of the nested compact sets in (1), and set \( p = L^{-1}p^{L} \).
\end{proof} 
In the rest of this section, \( \sigma \) is a solution of ABGCF on \( [0,T_{0}] \) with \( r_{0} \le \sigma_{t} \le R_{0} \), and \( L(0) = \mathrm{Id} \). By \cref{xiestimate},
\begin{equation}
    \int_{0}^{t}|\Lambda(s)|\,ds \le C\Big( t + t^{\frac{1-n\alpha}{1+n\alpha}} + 1 \Big), \label{Lambdaintegral}
\end{equation}
where \( C \) depends only on \( n \), \( \alpha \), \( Z_{1}(\sigma_{0}) \), \( r_{0} \) and \( R_{0} \).
The next lemma controls the geometry of \( X^{L} \) assuming bounds on the geometry of the original flow \( X \).
\begin{lemma}
    There exists a constant \( C = C(n,\alpha,Z_{1}(\sigma_{0}),r_{0},R_{0}) > 0 \) such that
    \[ C^{-1}e^{-Ct} \le \Psi(L(t),z) \le Ce^{Ct} \]
    and
    \[ B\big( 0, r_{0}C^{-1}e^{-Ct} \big) \subset L(t)\Omega_{t} \subset B\big( 0, R_{0}Ce^{Ct} \big) \]
    for all \( t \in [0,T_{0}] \).\label{PsiandSigmaL}
\end{lemma}
\begin{proof}
    By \cref{Lambdaintegral}, \( \|L(t)\| + \|L(t)^{-1}\| \le 2\exp\int_{0}^{t}|\Lambda| \le Ce^{Ct} \), using \( t^{\frac{1-n\alpha}{1+n\alpha}} \le 1+t \). The bounds on \( \Psi \) follow from the remark after \cref{supportL}. Since \( B(0,r_{0}) \subset \Omega_{t} \subset B(0,R_{0}) \), the image \( L(t)\Omega_{t} \) contains \( B(0,r_{0}/\|L(t)^{-1}\|) \) and is contained in \( B(0,R_{0}\|L(t)\|) \).
\end{proof}

For fixed \( \Lambda(t) \) (hence fixed \( L(t) \)), the \cref{supportL} is a parabolic equation which satisfies avoidance principle. In particular, we can use barriers to compare displacements. A lower displacement bound for \cref{supportL} can now be obtained following \cite[Theorem 8]{andrews2000motion}.
\begin{theorem}\label{displacement}
    There exist constants \( C, T' >0 \) depending continuously on \( n,\alpha,Z_{1}(\sigma_{0}),r_{0} \) and \( R_{0} \) such that
    \[ \sigma^{L}(z,t) - \sigma^{L}(z,0) \le -Ct^{\frac{1}{1-n \alpha}} \]
    for \( z \in \Sp^{n} \) and \( t \in [0,T'] \).
\end{theorem}
\begin{proof}
Let \( T=T(n,\alpha,Z_{1}(\sigma_{0}),r_{0},R_{0}) \) be such that \( e^{CT}\le 2 \) so that \( (2C)^{-1} \le \Psi \le 2C  \) for \( t \in [0,T] \) where \( C \) is the constant from \cref{PsiandSigmaL}. Let \( z \in \Sp^{n} \) and \( X^{L}(z,0) = \sigma^{L}(z)z + \nabla \sigma^{L}(z) \) be the point on the hypersurface with normal \( z \). The hypersurface \( \mathcal{M}^{L}_{0}= \mathcal{M}_{0} \) is enclosed in a hemispherical region of radius \( 2R_{0} \) centered at \( X^{L}(z) \) with the half-space with \( \{y \in \R^{n+1} : y \cdot z \le \sigma(z)\} \). Let \( \epsilon < 2 R_{0} \), we can find a sphere \( S_{\epsilon} \) with center at \( X^{L}(z) - \frac{4 R_{0}^{2}- \epsilon^{2}}{2 \epsilon}z \) and radius \( R_{\epsilon} = \frac{4 R_{0}^{2}+ \epsilon^{2}}{2 \epsilon} \) which encloses the hemispherical region and its support function in the direction \( z \) has values \( \sigma^{L}(z) + \epsilon \). The sphere \( S_{\epsilon} \)(t) evolving by
    \[ \frac{d}{dt}R(t)  = -(2C)^{-1}R(t)^{-n \alpha} \]
    with \( R(0) = R_{\epsilon} \) bounds \( \mathcal{M}_{t}^{L} \). Furthermore, \[ R(t) = \left( R_{\epsilon}^{1+n \alpha} - (1+n \alpha)(2C)^{-1}t \right)^{\frac{1}{1+n \alpha}} \le R_{\epsilon} - R_{\epsilon}^{-n \alpha}(2C)^{-1}t.\]
    By comparison principle, hypersurface \( \mathcal{M}^{L}_{t} \) is enclosed in \( S_{\epsilon}(t) \) and so
    \[ \sigma^{L}(z,t) \le \sigma^{L}(z,0) + \epsilon + R(t)-R_{\epsilon} \le \sigma^{L}(z,0) + \epsilon -  R_{\epsilon}^{-n \alpha}(2C)^{-1}t. \]
    It follows that, \begin{align*}
        \sigma^{L}(z,t) - \sigma^{L}(z,0) & \le \epsilon -  \left( \frac{2 \epsilon}{4 R_{0}^{2} + \epsilon^{2}} \right)^{n \alpha}(2C)^{-1}t \\
        & \le \epsilon - \left( \frac{\epsilon}{4 R_{0}^{2}} \right)^{n \alpha}(2C)^{-1}t \\
        & = \epsilon \left( 1- (4 R_{0}^{2})^{-n \alpha}\epsilon^{n \alpha-1}(2C)^{-1}t \right),
    \end{align*}
    where we used \( \epsilon < 2 R_{0} \) in the second last line. Choosing \( \epsilon = \left( 2^{-1-2 n \alpha}R_{0}^{-2n \alpha}(2C)^{-1}t \right)^{\frac{1}{1-n \alpha}} \) so that the term in bracket is negative one gives
    \[ \sigma^{L}(z,t) - \sigma^{L}(z,0) \le - \epsilon = -\left( 2^{-1-2 n \alpha}R_{0}^{-2 n \alpha}(2C)^{-1}t \right)^{\frac{1}{1-n \alpha}} \]
    for \( t \in [0,T') \) where \( T'=T'(n,\alpha,Z_{1}(\sigma_{0}),r_{0},R_{0}) \le T \) is chosen such that \( \epsilon < 2 R_{0} \) for \( t < T' \).
\end{proof}
Pairing up the displacement bound with Smoczyk type argument, we get the following lower speed bound.
\begin{theorem}\label{lowerspeed}
    There exist constants \( C , T' > 0 \) depending continuously on \( n,\alpha, Z_{1}(\sigma_{0}),r_{0} \) and \( R_{0} \) such that
    \[  Q_{n}[\sigma^{L}]^{-\alpha}(z,t) \ge C t^{ \frac{n \alpha}{1-n \alpha}} \]
    for every \( z \in \Sp^{n} \) and \( t \in (0,T'] \).
\end{theorem}
\begin{proof}
    The time evolution of speed \( Q_{n}^{-\alpha} \) under \cref{supportL} is given by 
    \begin{align*}
        \partial_{t}Q_{n}[\sigma^{L}]^{-\alpha} & = -\alpha Q_{n}^{-\alpha}(A^{-1})^{ij}A[\partial_{t}\sigma^{L}]_{ij} \\
        & = \mathcal{L}(\Psi Q_{n}^{-\alpha}) + \alpha \Psi H Q_{n}^{-2 \alpha}.
    \end{align*}
    Also, \( \mathcal{L}\sigma^{L} = n \alpha Q_{n}^{-\alpha} - \alpha H Q_{n}^{-\alpha}\sigma^{L} \). The evolution of the anisotropic factor is given by
    \begin{align*}
        |\partial_{t}\log \Psi| & = (1-\alpha(n+2)) \frac{|\left< \partial_{t}L^{T}z,L^{T}z \right>|}{|L^{T}z|^{2}} \\
        & \le (1-\alpha(n+2))|\xi(t)| \le |\xi(t)| \\
        & \le C(n,\alpha,Z_{1}(\sigma_{0}), r_{0},R_{0})\left(1 + t^{ -\frac{2n \alpha}{1+n \alpha}} + \frac{d}{dt}\log Z_{1}(\sigma)\right).
    \end{align*}
    Let \( \omega(t) = t (Z_{1}(\sigma_{t}) / Z_{1}(\sigma_{0}))^{\beta} \ge t \) and \( \Phi  = \sigma^{L}+ (2+n \alpha)\omega(t)\Psi Q_{n}[\sigma^{L}]^{-\alpha}  \)
    where \( \beta > 0 \) will be chosen later (note that we're taking entropy of the original flow, not of \cref{supportL}). The evolution of \( \Phi \) is given by
    \begin{align*}
        \partial_{t}\Phi & = -\Psi Q_{n}^{-\alpha} + (2+n \alpha)\left( \partial_{t}\omega(t)\Psi Q_{n}^{-\alpha} + \omega(t)\partial_{t}\Psi Q_{n}^{-\alpha}  + \omega(t)\Psi\mathcal{L}(\Psi Q_{n}^{-\alpha}) + \alpha \omega(t)\Psi^{2} H Q_{n}^{-2 \alpha}\right) \\
        & = \Psi\mathcal{L}\Phi + \alpha \Psi H Q_{n}^{-\alpha} \Phi + (2+n \alpha)\partial_{t}\omega \Psi Q_{n}^{-\alpha} - (1+n \alpha)\Psi Q_{n}^{-\alpha} + (2+n \alpha)\omega \partial_{t}\Psi Q_{n}^{-\alpha}\\
        & \ge \Psi \mathcal{L}\Phi + \alpha \Psi H Q_{n}^{-\alpha}\Phi +  (2+n \alpha)\left( (Z_{1}(t) / Z_{1}(0))^{\beta} + \omega \beta \frac{d}{dt}\log Z_{1} \right)\Psi Q_{n}^{-\alpha} - (1+n \alpha) \Psi Q_{n}^{-\alpha} \\
        & \quad - (2+n \alpha)\omega  C\left(1+t^{ \frac{-2n \alpha}{1+n \alpha}} + \frac{d}{dt}\log Z_{1}\right) \Psi Q_{n}^{-\alpha}.
    \end{align*}
    Choose \( \beta = C \) to cancel the entropy derivative terms and \( T' = T'(n,\alpha,Z_{1}(\sigma_{0}),r_{0},R_{0}) \) such that \newline \( (2+n \alpha)C(T'+T'^{ \frac{1-n \alpha}{1+n \alpha}}) \le 1 /2 \). Then for \( 0 \le t \le T' \), the evolution equation satisfies \begin{align*}
        \partial_{t}\Phi & \ge \Psi \mathcal{L}\Phi + \alpha \Psi H Q_{n}^{-\alpha}\Phi + (2+n \alpha)(Z_{1}(t) / Z_{1}(0))^{\beta}\Psi Q_{n}^{-\alpha} -(1+n \alpha)\Psi Q_{n}^{-\alpha}  \\
        & \quad - (2+n \alpha)  C(t + t^{ \frac{1-n \alpha}{1+n \alpha}}) (Z_{1}(t) / Z_{1}(0))^{\beta}\Psi Q_{n}^{-\alpha}  \\
        & \ge \Psi \mathcal{L}\Phi + \alpha \Psi H Q_{n}^{-\alpha}\Phi + (1+n \alpha)\Psi Q_{n}^{-\alpha} \left( (Z_{1}(t) / Z_{1}(0))^{\beta} - 1 \right)  \\
        & \quad + \left( 1 - (2+n \alpha)  C(t + t^{ \frac{1-n \alpha}{1+n \alpha}}) \right)(Z_{1}(t) / Z_{1}(0))^{\beta}\Psi Q_{n}^{-\alpha} \\
        & \ge \Psi \mathcal{L}\Phi + \alpha \Psi H Q_{n}^{-\alpha}\Phi + \frac{1}{2}(Z_{1}(t) / Z_{1}(0))^{\beta}\Psi Q_{n}^{-\alpha}.
    \end{align*} 
    Fix \( z' \in \Sp^{n} \). Since \( \Psi \) depends only on the normal and \( Q_{n}[\sigma^{L} + \langle p,z\rangle] = Q_{n}[\sigma^{L}] \), equation \cref{supportL} and the displacement bound of \cref{displacement} are invariant under translations, and so are the bounds on \( \Psi \) from \cref{PsiandSigmaL}. Translating the origin to \( X^{L}(z',0) \), we may therefore assume that \( \sigma^{L}(\cdot,0) \ge 0 \) and \( \sigma^{L}(z',0) = 0 \). Applying maximum principle, we conclude that \( \Phi(z,t)\ge 0 \) for \( 0\le t \le T' \) and hence  
    \[    Q_{n}[\sigma^{L}](z',t)^{-\alpha} \ge  - \frac{Z_{1}(\sigma_{0})^{\beta}\sigma^{L}(z',t)}{(2+n \alpha)t\Psi Z_{1}(\sigma_{t})^{\beta} } \ge  Ct^{ \frac{n \alpha}{1-n \alpha}}. \]
\end{proof}

Working in a time interval for which \( \Psi \) is bounded above and below, the lower speed bound for \( \sigma^{L} \) together with \cref{affineQ} immediately gives a lower speed bound for the original flow \( \sigma \). 
\begin{corollary}\label{lowerspeedoriginal}
    There exist constants \( C,T' >0  \) depending continuously on \( n, \alpha, Z_{1}(\sigma_{0}), r_{0} \) and \( R_{0} \) such that
    \[ Q_{n}[\sigma](z,t)^{-\alpha}  \ge C t^{ \frac{n \alpha}{1-n \alpha}} \]
    for \( z \in \Sp^{n} \) and \( t \in (0,T'] \).
\end{corollary}
Now we bound the individual principal curvatures for \( n > 1 \); for \( n = 1 \) the curvature is \( Q_{1}^{-1} \), so the speed bounds already bound it.
\begin{theorem}\label{curvaturebounds}
    Let \( n > 1 \) and suppose in addition that \( C_{1}^{-1} \le Q_{n} \le C_{1} \) on \( \Sp^{n}\times[0,T_{0}] \). Then there exist constants \( C_{2}, C_{3},C_{4}, C_{5} >0 \) depending continuously on \( n,\alpha,Z_{1}(\sigma_{0}),r_{0},R_{0} \) and \( C_{1} \) such that
    \[ \min \{C_{2},C_{3}t^{n-1}\} \le \kappa_{\min}(z,t) \le \kappa_{\max}(z,t) \le \max\{C_{4},C_{5}t^{-(n-1)^{2}}\} \]
    for every \( z \in \Sp^{n} \) and \( t \in (0,T_{0}] \).
\end{theorem}
\begin{proof}
    We analyze the evolution of \( A[\sigma]_{ij} \), whose eigenvalues are the principal radii. By \cref{evolution}, and since the terms in the second line are non-positive,
    \begin{align*}
    \partial_{t} A[\sigma]_{ij} & = \mathcal{L}A_{ij} + (n \alpha - 1)Q_n^{-\alpha}g_{ij}  - \alpha H Q_n^{-\alpha}A_{ij} \\
    & \quad - \alpha^{2}Q_n^{-2-\alpha}\nabla_{i}Q_n\nabla_{j}Q_n - \alpha Q_n^{-\alpha}(A^{-1})^{kp}(A^{-1})^{lq}\nabla_{i}A_{kl}\nabla_{j}A_{pq} \\
    & \quad + A_{ij}\xi + \frac{1}{2}\left< \nabla A_{ij}, \nabla \xi \right> + \frac{1}{2}g^{pq}A_{ip}\nabla_{q}\nabla_{j}\xi + \frac{1}{2}g^{pq}A_{jp}\nabla_{q}\nabla_{i}\xi \\
    & \le \mathcal{L}A_{ij} + (n \alpha - 1)Q_n^{-\alpha}g_{ij}  - \alpha H Q_n^{-\alpha}A_{ij}  + A_{ij}\xi + \frac{1}{2}\left< \nabla A_{ij}, \nabla \xi \right> \\
    & \quad + \frac{1}{2}g^{pq}A_{ip}\nabla_{q}\nabla_{j}\xi + \frac{1}{2}g^{pq}A_{jp}\nabla_{q}\nabla_{i}\xi.
\end{align*}
Let \( r_{\max}(t) = \sup_{z \in \Sp^{n}}\sup_{v \in T_{z}\Sp^{n}}A[\sigma]_{vv}(z) /g_{vv} \) be the maximum eigenvalue of \( A \) over \( \Sp^{n} \). At a point where the maximum is attained, we choose coordinates such that \( A \) and \( g \) are diagonal, \( v = \partial_{1} \) and \( \nabla A_{11} = 0 \), and we estimate \( |\xi| + |\nabla_{1}\nabla_{1}\xi| \le C'(n)|\Lambda(t)| \). Since \( Q_{n}^{-2\alpha} \le C_{1}^{2\alpha} \), \cref{xilemma} gives a constant \( C= C(n,\alpha,r_{0},R_{0},Z_{1}(\sigma_{0}),C_{1}) \) such that
\[ |\Lambda(t)| \le C \left( 1+ \frac{d}{dt}\log Z_{1}\right). \]
Thus,
\begin{align*}
    \frac{d}{dt}r_{\max} & \le (n \alpha-1) Q_{n}^{-\alpha} - \alpha H Q_{n}^{-\alpha}r_{\max} + C'(n)|\Lambda(t)|r_{\max} \\
    & \le  (n \alpha-1) Q_{n}^{-\alpha} - \alpha H Q_{n}^{-\alpha}r_{\max} + C'C\left(1 + \frac{d}{dt}\log Z_{1}\right)r_{\max}.
\end{align*}
We also have \( Q_{n} = \prod r_{i} \ge r_{\min}^{n-1}r_{\max} \), so \( H = \sum r_{i}^{-1} \ge r_{\min}^{-1} \ge Q_{n}^{ -\frac{1}{n-1}}r_{\max}^{\frac{1}{n-1}} \ge C_{1}^{-\frac{1}{n-1}}r_{\max}^{\frac{1}{n-1}}  \), and \( n \alpha - 1 < 0 \), so
\[ \frac{d}{dt}r_{\max} \le - \alpha C_{1}^{-\alpha - \frac{1}{n-1}}r_{\max}^{\frac{n}{n-1}} + C'C\left(1 + \frac{d}{dt}\log Z_{1}\right)r_{\max}.
\]
Let \( \psi = r_{\max}Z_{1}^{-C'C} \). Then, using \( Z_{1} \ge Z_{1}(\sigma_{0}) \),
\begin{align*}
    \frac{d}{dt}\psi & \le -\alpha C_{1}^{-\alpha - \frac{1}{n-1}}\psi^{\frac{n}{n-1}}Z_{1}^{\frac{C'C}{n-1}} + C'C \psi \\
    & \le -C''\psi^{\frac{n}{n-1}}+ C'C \psi.
\end{align*}
By the maximum principle,
\[ \psi \le \max \left\{(2C'C C''^{-1})^{n-1}, \left( \frac{C''t}{2(n-1)} \right)^{-(n-1)}\right\}, \] and thus, by \cref{Z1upperbound},
\[ r_{\max} \le \max\{C_{2}^{-1},C_{3}^{-1}t^{-(n-1)}\}. \]
This gives the lower bound on \( \kappa_{\min} = r_{\max}^{-1} \). The upper bound on \( r_{\max} \) implies a lower bound on \( r_{\min} \) from \( r_{\max}^{n-1}r_{\min} \ge Q_{n} \ge C_{1}^{-1} \), which gives the upper bound on \( \kappa_{\max} = r_{\min}^{-1} \).
\end{proof}

\section{Sequential convergence}

With the curvature estimates from previous section, we can prove that the flow converges to a point for a sequence of times. First, we need estimates on a window which will be used to prove convergence in the time dependent anisotropic formulation from \cref{sec:lowerspeed}. 

Throughout this section, \( X \) is a solution of ABGCF with maximal time \( T \). The letter \( C \) will denote constants depending only on \( n, \alpha \) and \( Z_1(\Omega_0) \). Let \( K_{Z} = \log\big( C_{\mathrm{iso}}/Z_{1}(\Omega_{0}) \big) \), where \( C_{\mathrm{iso}} \) is the constant of \cref{Z1upperbound}. Then \( \int_{t_{1}}^{t_{2}}\frac{d}{dt}\log Z_{1}\,dt \le K_{Z} \) for all \( t_{1} \le t_{2} \), by \cref{Monotonicity}. The constants in \cref{sec:upperspeed,sec:lowerspeed} depend on \( Z_{1}(\sigma_{0}) \) only through a lower bound for it, so when we apply those results from an initial time \( t_{1} \) we may use \( Z_{1}(\Omega_{0}) \le Z_{1}(\Omega_{t_{1}}) \) in place of \( Z_{1}(\Omega_{t_{1}}) \).

\begin{proposition}[Estimates on a window]\label{window}
    There exists \( h \in (0,1] \) and \( C \) with the following property. Let \( t_1 \in [0,T) \) with \( V(\Omega_{t_{1}}) = V(B^{n+1}) \), let \( t_{2} = \min \{t_{1}+h,T\} \) and let \( L' = -L \Lambda \) on \( [t_{1},t_{2}) \) with \( L(t_{1}) = \mathrm{Id} \). Then \begin{enumerate}
    \item \( ||L(t)||, ||L(t)^{-1}|| \le C \) and \( \int_{t_{1}}^{t_{2}}|\Lambda| \le C \) for \( t \in [t_{1},t_{2}) \),
    \item there is a curve \( q \) with \( q' = \Lambda q \) such that \( B(q(t),C^{-1}) \subset \Omega_{t} \subset B(q(t),C) \) for \( t \in [t_{1},t_{2}) \),
    \item for \( \delta \in (0, h /2] \), there exists a \( C_{\delta} \) such that 
    \[ C_{\delta}^{-1} \le Q_{n}[\sigma_{t}] \le C_{\delta},\quad \text{ and } \quad C_{\delta}^{-1}g \le A[\sigma_{t}] \le C_{\delta}g  \text{ on } t \in [t_{1}+ \delta, t_{2}).\]
\end{enumerate}
\end{proposition}
\begin{proof}
    (1) The isoperimetric ratio is controlled by \cref{Volumeandgeometry} and \cref{isoperimetricbounded}, so for sufficiently small \( h \), the convex body \( \Omega_{t} \) contains a ball of radius \( r_{*} \) and is contained in a ball of radius \( R_{*} \) where \( r_{*} \) and \( R_{*} \) only depend on \( n, \alpha \) and \( Z_{1}(\Omega_{0}) \). 

    To establish (1) and (2), we will do a continuity argument. Suppose that there exists a time interval \( [t_{1},t'] \subset [t_{1},t_{2}) \) and a curve \( q \) with \( q' = \Lambda q \) and \( B(q(t),r_{0}) \subset \Omega_{t} \subset B(q(t),R_{0}) \) for \( t \in [t_{1},t'] \). By  \cref{translationlemma}(3), \( X - q \) is a solution of ABGCF with same \( \Lambda \), and moreover satisfies the hypothesis of \cref{upperspeed} and on \( [t_{1},t'] \), so using \cref{xiestimate}
    \begin{equation}
        \int_{t_{1}}^{t'}|\Lambda| \le C(r_{0},R_{0})\left(h + h^{\frac{1-n\alpha}{1+n\alpha}} \right) + C_{2}K_{Z}.\label{Lambdaint2}
    \end{equation}
    Let \( M = 4e^{C_{2}K_{Z}} \) and decrease \( h \) so that the first term in \cref{Lambdaint2} with \( r_{0} = r_{*} / M^{2} \) and \( R_{0} = 2R_{*} \) is at most \( \log 2 \). Define 
    \[ t^{*} = \sup \{t \in [t_{1},t_{2}) : ||L(s)|| \le M \text{ and } ||L(s)^{-1}|| \le M \text{ for }s \in [t_{1},t]\} > t_{1}. \]
    We claim that \( t^{*} = t_{2} \). Fix \( t' < t^{*} \) and a ball \( B(x',r_{*}) \subset \Omega_{t'} \), then for \( s \in [t_{1},t'] \), the bodies \( L(s)\Omega(s) \) are nested by \cref{nesting} so that 
    \[ B(L(t')x', r_{*} /M) \subset L(t')B(x',r_{*}) \subset L(t')\Omega_{t'} \subset L(s)\Omega_{s}. \]
    It follows that 
    \[ B(q_{t'}(s), r_{*} / M^{2}) \subset \Omega_{s}, \quad \text{ where } \quad  q_{t'}(s) = L(s)^{-1}L(t')x' \]
    and also the time derivative of the curve \( q_{t'} \) satisfies \( q_{t'}' = \Lambda q_{t'} \). Now \( q_{t'}(s) \in \Omega_{s} \), and \( \Omega_{s} \subset B(q_{t'}(s),2R_{*}) \) so by the earlier estimate, we have 
    \[ \int_{t_{1}}^{t'} |\Lambda| \le \log 2 + C_{2}K_{Z} = \log (M /2), \]
    and thus by the setup in \cref{supportL}, \( ||L(t)^{\pm}|| \le \exp \int_{t_{1}}^{t'} |\Lambda| \le M /2 \) on \( [t_{1},t^{*}) \). By the definition of \( t^{*} \), we must have \( t^{*} = t_{2} \) which establishes (1).

    For (2), let \( t' \to t_{2} \) and suppose \( q_{t'}(t_{1}) \in \Omega_{t_{1}} \) converges to some \( q_{1} \). The curve \( q(t) = L(t)^{-1}q_{1} \) with \( C = \max\{M^{2} /r_{*},2R_{*}\} \) satisfies the property.

    For (3), we apply curvature estimates from \cref{sec:lowerspeed,sec:upperspeed} to \( X - q \) starting at time \( t_{1} + \delta / 2 \) after decreasing \( h \) enough such that \( h \le T' \), where \( T' \) is the time from \cref{lowerspeedoriginal} with \( r_{0} = r_{*} / M^{2} \) and \( R_{0} = 2R_{*} \). 

\end{proof}

The above proposition allows us to establish higher regularity on a window. We will work with the time-dependent anisotropic formulation from \cref{sec:lowerspeed}. 

\begin{proposition}\label{higherregularity}
    Under the assumptions of \cref{window}, for every \( k \) and \( \delta \in (0, h /2] \) there exists a \( C_{\delta,k} \) such that \( ||\sigma_{t} - \left< q(t),z \right>||_{C^{k}(\Sp^{n})}\le C_{\delta,k} \) for all \( t \in [t_{1} + 2 \delta,t_{2}) \) for \( q(t) \in \bar{\Omega}_{t} \) from (2) of \cref{window}.
\end{proposition}
\begin{proof}
    Let \( \sigma^{L} \) be the support function of \( L(t)(\Omega_{t} - q(t)) \). Then from \cref{supportL}, \( \sigma^{L} \) satisfies the parabolic PDE
    \[ \partial_{t}\sigma^{L} = - \Psi Q_{n}[\sigma^{L}]^{-\alpha}, \quad \text{ with } \quad \Psi(z,t) = |L(t)^{T}z|^{1-\alpha(n+2)}. \]
    By \cref{window}, \( \|L(t)^{\pm 1}\| \le C \) on \( [t_{1},t_{2}) \), so \( C^{-1} \le \Psi \le C \) and \( \|\Psi(\cdot,t)\|_{C^{k}} \le C_{k} \) for every \( k \). Moreover, by \cref{window}(3), \( A[\sigma^{L}] \) is pinched on \( [t_{1}+\delta,t_{2}) \).

    The speed \( f = Q_{n}[\sigma^{L}]^{-\alpha} \) satisfies
    \[ \partial_{t}f = \alpha f (A^{-1})^{ij}A[\Psi f]_{ij} = a^{ij}\nabla_{i}\nabla_{j}f + b^{i}\nabla_{i}f + cf, \]
    with \( a^{ij} = \alpha f\Psi(A^{-1})^{ij} \), \( b^{i} = 2\alpha f(A^{-1})^{ij}\nabla_{j}\Psi \) and \( c = \alpha f(A^{-1})^{ij}A[\Psi]_{ij} \). On \( [t_{1}+\delta,t_{2}) \) these coefficients are bounded and measurable and \( (a^{ij}) \) is uniformly elliptic, so by the Krylov-Safonov theorem \cite{krylov1980property}, \( f \) is H\"older continuous in space and time. For each fixed \( t \), \( \sigma^{L}(\cdot,t) \) solves the uniformly elliptic equation \( (\det A[\sigma^{L}])^{\frac{1}{n}} = f^{-\frac{1}{n\alpha}} \), whose operator is concave, so Caffarelli's interior estimates \cite{caffarelli1989interior} give \( \|\sigma^{L}(\cdot,t)\|_{C^{2,\beta}} \le C \). Hence, \( a^{ij}, b^{i}, c \) are H\"older continuous in space, uniformly in \( t \). The interior Schauder estimates for parabolic equations whose coefficients are H\"older continuous in space and only measurable in time give \( ||f(\cdot,t)||_{C^{2,\beta}} \le C \), and elliptic Schauder theory applied to the equation for \( \sigma^{L}(\cdot,t) \) then gives \( ||\sigma^{L}(\cdot,t)||_{C^{4,\beta}} \le C \). Differentiating the equation for \( f \) in space and repeating the argument on slightly smaller time intervals gives \( ||\sigma^{L}(\cdot,t)||_{C^{k}} \le C_{\delta,k} \) for every \( k \) and \( t \in [t_{1}+2\delta,t_{2}) \). Since \( ||L(t)^{\pm 1}|| \le C \), these bounds transfer to \( \sigma_{t} - \left< q(t),z \right> \), which is the support function of \( \Omega_{t} - q(t) \).
\end{proof}

Now, we're ready to prove that ABGCF exists for maximal time and converges to a point. 

\begin{theorem}\label{extinction}
    The maximal time \( T \) is finite, and \( \text{diam}(\Omega_{t}) \to 0 \) as \( t \to T \). With \( p(t) \) as in \cref{nesting}(3), \( \Omega_{t} - p(t) \) converges to the origin.
\end{theorem}
\begin{proof}
    By \cref{Volumeandgeometry}, \( T < \infty \). Suppose \( V(\Omega_{t}) \ge a >0 \) for all \( t < T \). As the isoperimetric ratio is controlled throughout the flow, we can find a uniform window \( h_{0} \) in the sense of \cref{window} for which the curvature estimates remain valid. But then applying the estimate at time \( T - h_{0} / 2 \) implies that \( ||L(t)|| \) remains bounded and the bodies \( \Omega_{t} \subset L(t)^{-1}\Omega_{0}\) stay uniformly convex so \( \sigma_{t} \) is bounded in every \( C^{k} \). Moreover, \( ||\partial_{t}\sigma||_{C^{k}} \le C_{k}(1+|\Lambda|) < \infty \) so \( \sigma_{t} \) converges in \( C^{\infty} \) to a uniformly convex, \( Z_{1} \)-isotropic limit. By \cref{shorttime}, the solution can be extended beyond \( T \), and we arrive at a contradiction. Hence, \( V(\Omega_{t}) \to 0 \) and by \cref{isoperimetricbounded}, \( \text{diam}(\Omega_{t}) \to 0 \).
\end{proof}

We can also prove convergence to a self-similar solution for a sequence of time. Let \( p \) be the curve from \cref{nesting}(3) on the maximal time interval \( [0,T) \), then \( X-p \) is also a solution of ABGCF, converging to the origin. We rescale to fixed volume as 
\begin{equation}
    \rho(t) = \left( \frac{V(B^{n+1})}{V(\Omega_{t})} \right)^{\frac{1}{n+1}}, \qquad \tilde{X}_{t} = \rho(t)X_{t}. \label{rescaled}
\end{equation}

By \cref{Volumeandgeometry} and \cref{extinction},  \begin{equation}
    C^{-1}(T-t)^{-1} \le \rho(t)^{1+n \alpha} \le C(T-t)^{-1}. \label{rhoestimate}
\end{equation}

\begin{theorem}
    There exists a sequence of time \( t_{i} \to T \) and a uniformly convex self-similar solution (which is in \( Z_{1} \) isotropic position by criticality with respect to \( Z_{1} \) entropy) with support function \( \omega \) such that \( \tilde{\sigma}_{t_{i}} \to \omega \) in \( C^{\infty} \). \label{seqconvergence}
\end{theorem}

\begin{proof}
    
    Fix \( t_{1} < T \) and let \( Y(\cdot,s) = \rho(t_{1})X\big( \cdot,t_{1} + \rho(t_{1})^{-(1+n\alpha)}s \big) \). By \cref{Volumeandgeometry}, after decreasing \( h \) we have \( V(\Omega_{t}) \ge \frac{1}{2}V(\Omega_{t_{1}}) \) for \( t \le t_{1} + h\rho(t_{1})^{-(1+n\alpha)} \), so this time is less than \( T \) and \( \rho(t)/\rho(t_{1}) \in [1,2^{\frac{1}{n+1}}] \). \Cref{window,higherregularity}, applied to \( Y \) with \( \delta = h/4 \), bound \( \tilde{X}(\cdot ,t ) = \rho(t) / \rho(t_{1}) Y( \cdot, s) \) in every \( C^{k} \) and give \( A[\tilde\sigma_{t}] \ge C^{-1}g \) for \( t \in \big[ t_{1} + \frac{h}{2}\rho(t_{1})^{-(1+n\alpha)},\ t_{1} + h\rho(t_{1})^{-(1+n\alpha)} \big) \). As \( t_{1} \to T \), the left endpoint of this interval moves continuously to \( T \), so these intervals cover \( [t',T) \) for some \( t' < T \).

    By a standard argument using Arzel\`a-Ascoli theorem, there exists a convergent subsequence in \( C^{\infty} \). The final convex body is a critical point of \( Z_{1} \) entropy, and hence a self-similar solution of \( \alpha \)-GCF. The details of the argument can be found in \cite[Section 9]{andrews1997monotone}. Since \( Z_{1} \) is invariant under dilations and translations, continuous under \( C^{2} \) convergence and non-decreasing by \cref{Monotonicity}, we have \( Z_{1}(\omega) = \lim_{i \to \infty} Z_{1}(\mathcal{M}_{t_{i}}) \ge Z_{1}(\mathcal{M}_{0}) \). 

\end{proof}

\section{Construction of symmetric self-similar solutions} \label{construction}

We prove the main result in this section. The first author showed the existence of self-similar solution possessing an irreducible group of symmetry in \cite{andrews2000motion}. ABGCF allows us to construct self-similar solutions with weaker symmetry assumptions. In particular, we can establish the existence of solutions with reducible symmetry groups like axial symmetry and discrete subgroups like dihedral symmetry groups.

\begin{theorem}
    Let \( G \subset O(n+1) \) be a proper, closed subgroup for which there exists a spherical harmonic \( \phi \) with degree \( d \ge 3 \) such that \( \phi \) is \( G \)-symmetric. Then for any \( 0 < \alpha < \frac{1}{(d-1)(n+d)} \) there exists a non-spherical, compact, locally uniformly convex self-similar solution to the \( \alpha \)-GCF with symmetry group containing \( G \).\label{Ginvariantsoln}
\end{theorem}
\begin{proof}
    Let \( m_{d} = (d-1)(n+d) \) and \( \Omega_{\epsilon} \) be a fixed volume variation of the ball with support function \( \sigma_{\epsilon}\) such that \( \sigma_{0} = 1 \), \( \dot\sigma = \phi \) and \( \Omega_{\epsilon} \) is \( G \)-invariant. Note that \( \Delta \phi = -d(d+n-1)\phi \), so \( (\Delta + n)\phi = - m_{d}\phi \) and the condition on \( \alpha \) is equivalent to \( 1- \alpha m_{d} > 0 \).

    \textit{Step 1:} Let \( T_{\epsilon} \in \operatorname{SL}(n+1) \) from \cref{isotropicZ} which is symmetric and minimizes \( Z_{1}(T_{\epsilon}\Omega_{\epsilon}) \) so that \( T_{\epsilon}\Omega_{\epsilon} \) is in \( Z_{1} \)-isotropic position. We claim that \( T_{\epsilon}\Omega_{\epsilon} \) is \( G \)-invariant as well. Let \( g \in G \subset O(n+1)\), then \( gTg^{-1} \) is symmetric as well and \( g^{-1}\Omega_{\epsilon} = \Omega_{\epsilon} \) by hypothesis, so we have 
    \[ Z_{1}(gT_{\epsilon}g^{-1}\Omega_{\epsilon}) = Z_{1}(gT_{\epsilon}\Omega_{\epsilon}) = Z_{1}(T_{\epsilon}\Omega_{\epsilon}), \]
    where the last equality follows from rotational invariance of \( Z_{1} \) entropy. By the uniqueness in \cref{isotropicZ}, we must have 
    \[ gT_{\epsilon}g^{-1} = T_{\epsilon}, \quad \forall g \in G. \]
    Hence, for any \( g \in G \), we have \[ gT_{\epsilon}\Omega_{\epsilon} = T_{\epsilon}g\Omega_{\epsilon} = T_{\epsilon}\Omega_{\epsilon}. \] 
    Thus, \( T_{\epsilon}\Omega_{\epsilon} \) is \( G \)-invariant.
    
        \textit{Step 2:} Let \( \omega_{\epsilon} \) denote the support function of \( T_{\epsilon}\Omega_{\epsilon} \). Since the ball is in \( Z_{1} \)-isotropic position, \( T_{0} = \mathrm{Id} \), and since \( D^{2}\tilde{Z}_{1} \) is non-degenerate by \cref{HessianZpq}, the implicit function theorem shows that \( T_{\epsilon} \) depends smoothly on \( \epsilon \). Hence \( \omega_{0} = 1 \) and, by \cref{DAsigma}, \( \psi \coloneqq \dot\omega_{0} = \phi + \xi \), where \( \xi = \langle \dot{T}_{0}z,z \rangle \) is a (possibly zero) quadratic spherical harmonic. Note that \( (\Delta+n)\xi = -(n+2)\xi \) and \( \int_{\Sp^{n}}\phi\xi = 0 \).
    Now, the volume is preserved so
    \[ 0 = \frac{d^{2}}{d \epsilon^{2}}\bigg|_{\epsilon = 0}V(T_{\epsilon}\Omega_{\epsilon}) = \frac{d}{d \epsilon}\bigg|_{\epsilon = 0}\int_{\Sp^{n}}\dot \omega_{\epsilon}Q_{n}[\omega_{\epsilon}] = \int_{\Sp^{n}} \ddot \omega_{0} + \psi(\Delta + n)\psi. \]
    The first variation of \( \int Q_{n}^{1-\alpha} \) at the sphere is \( n(1-\alpha)\int\psi = 0 \). For the second variation,
    \begin{align*}
        \frac{d^{2}}{d \epsilon^{2}}\bigg|_{\epsilon = 0}\int_{\Sp^{n}}Q_{n}[\omega_{\epsilon}]^{1-\alpha} & = (1-\alpha)\frac{d}{d \epsilon}\bigg|_{\epsilon = 0}\int_{\Sp^{n}}nQ_{n}[\omega_{\epsilon}]^{-\alpha}Q[\dot \omega_{\epsilon}, \omega_{\epsilon},\dots, \omega_{\epsilon}]\\
        & = (1-\alpha) \int_{\Sp^{n}} -\alpha ((\Delta + n)\psi)^{2} + n\ddot \omega_{0} + (n-1)\psi (\Delta + n)\psi\\
        & = (1-\alpha) \int_{\Sp^{n}} -\alpha ((\Delta + n)\psi)^{2} - \psi (\Delta + n)\psi\\
        & = (1-\alpha)\left( m_{d}(1-\alpha m_{d})\int_{\Sp^{n}}\phi^{2} + (n+2)(1-\alpha(n+2))\int_{\Sp^{n}}\xi^{2} \right) \\
        & \ge (1-\alpha)m_{d}(1-\alpha m_{d})\int_{\Sp^{n}}\phi^{2},
    \end{align*}
    where we used \( \alpha < \frac{1}{n+2} \) in the last line. Since the volume is fixed, Taylor expansion gives
    \[ \log Z_{1}(\omega_{\epsilon}) \ge \left( \frac{1}{2}(n+1)(1-\alpha)m_{d}(1- \alpha m_{d})\fint_{\Sp^{n}}\phi^{2} \right)\epsilon^{2} + O(\epsilon^{3}). \]

    \textit{Step 3:} By Step 2, \( Z_{1}(T_{\epsilon}\Omega_{\epsilon}) > 1 = Z_{1}(B^{n+1}) \) for small \( \epsilon \neq 0 \). We run the ABGCF starting at the \( Z_{1} \)-isotropic hypersurface \( \mathcal{M}_{0} = \partial(T_{\epsilon}\Omega_{\epsilon}) \). For \( g \in G \), \( gX \) is a solution of ABGCF with the same initial hypersurface by \cref{translationlemma}(1), so \( gX = X \) by the uniqueness in \cref{shorttime}, and every \( \mathcal{M}_{t} \) is \( G \)-invariant. By \cref{seqconvergence}, \( \rho(t_{i})(\mathcal{M}_{t_{i}} - p(t_{i})) \to \mathcal{M}_{\infty} \), a self-similar solution with \( Z_{1}(\mathcal{M}_{\infty}) \ge Z_{1}(\mathcal{M}_{0}) > 1 \), so \( \mathcal{M}_{\infty} \) is not a sphere. Since \( p(t_{i}) \in \overline{\Omega}_{t_{i}} \) and \( \rho(t)\operatorname{diam}(\Omega_{t}) \le C \) by \cref{isoperimetricbounded}, after passing to a subsequence \( \rho(t_{i})p(t_{i}) \to p_{\infty} \), and then \( \rho(t_{i})\mathcal{M}_{t_{i}} \to \mathcal{M}_{\infty} + p_{\infty} \). A limit of \( G \)-invariant hypersurfaces is \( G \)-invariant, so \( \mathcal{M}_{\infty} + p_{\infty} \) is a non-spherical self-similar solution whose symmetry group contains \( G \).
\end{proof}
    \cref{maintheorem} now follows from the following algebraic fact about spherical harmonics.
    \begin{proposition}
        Let \( G \subset O(n+1) \) be a proper, closed subgroup that acts non-transitively on \( \Sp^{n} \), then there exists a spherical harmonic \( \phi \) with degree \( d \ge 3 \) such that \( \phi \) is \( G \)-symmetric. \label{harmonics}
    \end{proposition}
    \begin{proof}
    The group \( G \) acts non-transitively on \( \Sp^{n} \) via isometries, so the metric space \( \Sp^{n} / G \) is connected and has at least two distinct points and thus uncountably many points. By Urysohn's lemma, there exists a non-constant continuous function \( f : \Sp^{n} / G \to \R \). The pullback of \( f \) to \( \Sp^{n} \) is a non-constant, \( G \)-symmetric continuous function. Decompose \( f = \sum_{d \ge 0} f_{d} \) in \( L^{2}(\Sp^{n}) \), where \( f_{d} \) is a spherical harmonic of degree \( d \). The orthogonal projection onto the degree-\( d \) harmonics commutes with the action of \( O(n+1) \), so each \( f_{d} \) is \( G \)-symmetric. If every \( G \)-symmetric continuous function had \( f_{d} = 0 \) for all \( d \ge 3 \), then \( f \mapsto (f_{0},f_{1},f_{2}) \) would embed \( C(\Sp^{n}/G) \) into a finite-dimensional space. But \( C(\Sp^{n}/G) \) is infinite-dimensional, since \( \Sp^{n}/G \) is an infinite compact metric space. Hence some \( G \)-symmetric \( f \) has \( f_{d} \neq 0 \) for some \( d \ge 3 \), and \( \phi = f_{d} \) is the required harmonic.
    \end{proof}
    
    \textit{Proof of \cref{maintheorem}.} Combine \cref{Ginvariantsoln} and \cref{harmonics}.
\begin{corollary}\label{axial}
    There exists non-spherical, compact, locally uniformly convex self-similar solutions to the \( \alpha \)-GCF with axial symmetry for any \( 0< \alpha < \frac{1}{2(n+3)} \). 
\end{corollary}
\begin{proof}
    Let \( G \simeq O(n) \) be the group of orthogonal matrices fixing \( e_{n+1} \). There exists a spherical harmonic of degree \( 3 \), which is a cubic polynomial in \( z^{n+1} = \left< z, e_{n+1} \right> \) and is \( G \)-invariant. \cref{Ginvariantsoln} then gives the desired result. 
\end{proof}

\bibliographystyle{acm}
\bibliography{ref}
\end{document}